\documentclass{article}
\usepackage[latin1]{inputenc}
\usepackage[T1]{fontenc}
\usepackage{amstext}
\usepackage{amsmath, amssymb}
\usepackage[english]{babel}
\usepackage{amsfonts}

\newcommand{\R}{\mathbb{R}}
\newcommand{\C}{\mathbb{C}}
\newcommand{\N}{\mathbb{N}}

\newcommand\kD{\mathcal{D}}

\newcommand\kC{\mathcal{C}}
\newcommand\kR{\mathcal{R}}
\newcommand\kO{\mathcal{O}}

\newcommand\kF{\mathcal{F}}
\newcommand\kH{\mathcal{H}}
\newcommand\kI{\mathcal{I}}
\newcommand\kS{\mathcal{S}}

\newcommand\kL{\mathcal{L}}
\newcommand\ka{\mathfrak{a}}
\newcommand\kh{\mathfrak{h}}
\newtheorem {defin} {Definition} [section]
\newtheorem {lem} {Lemma} [section]
\newtheorem {prop} {Proposition} [section]

\newtheorem {theo} {Theorem} [section]
\newtheorem {cor} {Corollary} [section]
\newtheorem {rem} {Remark} [section]

\newcommand\la{\lambda}

\title{Contributions to the hypergeometric function theory of Heckman and Opdam:
sharp estimates, Schwartz space, heat kernel}
\author{Bruno Schapira\footnote{partially supported by the European Commission (IHP Network HARP $2002-2006$).}}

\begin{document}
\maketitle

\begin{center}
\it{Université d'Orléans,\\ Fédération Denis Poisson, Laboratoire MAPMO \\
B.P. 6759, 45067 Orléans cedex 2, France .}
\end{center}
\begin{center}
\it{Université Pierre et Marie Curie,\\
Laboratoire de Probabilités et Modèles Aléatoires, \\
4 place Jussieu, \\
F-75252 Paris cedex 05, France.}
\end{center}

\vspace*{0.8cm}

\begin{abstract}
Under the assumption of positive multiplicity, we obtain basic
estimates of the hypergeometric functions $F_\la$ and $G_\la$ of
Heckman and Opdam, and sharp estimates of the particular functions
$F_0$ and $G_0$. Next we prove the Paley-Wiener theorem for the
Schwartz class, solve the heat equation and estimate the heat
kernel. \end{abstract}

\bigskip
\noindent{\bf Key Words:} Differential-difference equations,
hypergeometric functions, root systems, Schwartz space, heat
kernel.

\bigskip
\noindent{\bf A.M.S. Classification.} {\it Primary:} 33C67,
33D67,
33E30, 42B10, 58J35. \\
{\it Secondary:} 35K05, 42A90, 43A32, 47D07, 58J65.

\bigskip
\noindent{\bf e-mail. } bruno.schapira@univ-orleans.fr

\vspace*{0.8cm}

\section{Introduction}
Classical harmonic analysis on $\R^n$ has now been extended to other spaces.
For instance Harish-Chandra has considered the case of semi-simple Lie groups.
Then he was followed by Helgason, who studied the Riemannian symmetric
spaces of noncompact type, which are Riemannian spaces of negative
curvature. In particular, Harish-Chandra introduced and studied
the spherical functions, which play the role of the exponentials
in these spaces. A more general setting, in the flat case, has
appeared two or three decades ago, with the theory of Dunkl
operators. It gives a vast generalization of the exponential
functions, and of the Fourier transform on $\R^n$. But it gives
also a generalization of the harmonic analysis on tangent spaces
of symmetric spaces. The natural counterpart of the Dunkl theory
in the negatively curved setting is the theory of Heckman and
Opdam. This theory has known a deep evolution with the discovery
of the Cherednik operators \cite{C}, the analogues of the Dunkl
operators in the flat case.
\newline
Heckman and Opdam \cite{HO}, \cite {HS}, \cite{O} have developed
their theory in the last two decades. They have first introduced a
new family of functions $F_\la$ on $\R^n$, which like in the Dunkl
theory are associated to root systems and a parameter, the
multiplicity function. They can be defined essentially as
eigenfunctions of certain differential operators. When the
multiplicity function, takes particular values, then these
operators coincide with the radial part of the $G$-invariant
differential operators on the symmetric spaces of noncompact
type $G/K$. Thus the restrictions to a Cartan subspace $\ka$ of
the spherical functions are particular functions $F_\la$. In this
way the theory of Heckman and Opdam is also a generalization of
the harmonic analysis on the symmetric spaces $G/K$. However all the
techniques used by Harish-Chandra can not always be transposed (at
least not trivially) in this new theory, because there are not anymore
underlying Lie groups. The main tools used in the harmonic
analysis on the symmetric spaces are in the one part an integral
formula of the spherical functions, and in another part a
development in series of these spherical functions. Heckman and
Opdam have shown that their functions $F_\la$ have a development
in series of the type Harish-Chandra, but there is not (at least
not yet) an integral formula, for general root systems. However
this gap has been compensated by two main discoveries. First the
discovery of the differential-difference operators by Cherednik
\cite{C}, and then the discovery by Opdam of a new type of
functions, the functions $G_\la$ \cite{O}, for which the calculus
and estimates can be more easily performed. These functions are
eigenfunctions of the Cherednik operators. However until recently
the only asymptotic result was essentially the fact that the
functions $F_\la$ and $G_\la$ were bounded \cite{O}. Delorme has
obtained a much better estimate, even in the more complicated case
of a negative multiplicity \cite{D}, but it requires involved
materials and techniques.
\newline In this paper we give sharp estimates of the functions $F_\la$, $G_\la$
and their derivatives, in an elementary way. Our method is only
based on the study of the system of differential and difference
equations satisfied by the functions $G_\la$, improving by the way
what had already done De Jeu \cite{J} and Opdam \cite{O} for
bounding their functions. We also give a global estimate of the
particular functions $F_0$ and $G_0$. It generalizes some results
in the noncompact symmetric spaces \cite{A}, \cite{ABJ}. Then we
deduce from these estimates and from a general method of Anker
\cite{A2} the inversion formula on the Schwartz space. Finally we
solve the Heat equation and we give some estimates of the heat
kernel.
\newline
\newline
\noindent \textbf{Acknowledgments:} This work is part of my PhD.
It is a great pleasure to thank my advisors Jean-Philippe Anker
and Philippe Bougerol for their help and advices.

\section{Preliminaries}

Let $\ka$ be a Euclidean vector space of dimension $n$, equipped
with an inner product $(\cdot ,\cdot )$. Let $\kh=\ka \otimes_\R
\C$ be the complexification of $\ka$. The notation $\Re$ and $\Im$
denote the real and imaginary part respectively, of an element in
$\kh$ or possibly in $\C$. Let $\kR \subset \ka$ be an integral
root system. We choose a subset of positive roots $\kR^+$. We
denote by $\kR_0^+$ the set of positive indivisible roots, by
$\Pi$ the set of simple roots, and by $Q^+$ the positive lattice
generated by $\kR^+$. Let $\alpha^\vee=\frac{2}{|\alpha|^2}\alpha$
be the coroot associated to a root $\alpha$ and let
$$r_\alpha(x)=x-(\alpha^\vee,x)\alpha,$$
be the corresponding orthogonal reflection. We denote by $W$ the
Weyl group associated to $\kR$, i.e. the group generated by the
$r_\alpha$'s. If $C$ is a subset of $\ka$, we call
\emph{symmetric} of $C$ any image of $C$ under the action of $W$.
Let $k\ :\ \kR \rightarrow [0,+\infty)$ be a multiplicity
function, which by definition is $W$-invariant. In the sequel we
may actually forget about the roots $\alpha$ with $k_\alpha =0$
and restrict ourself to the root subsystem where $k$ is strictly
positive.
\newline
Let
$$\ka_+ = \{x \mid \forall \alpha \in \kR^+,\ (\alpha,x)>0\},$$
be the positive Weyl chamber. We denote by $\overline{\ka_+}$ its
closure, and by $\partial \ka_+$ its boundary. Let also
$\ka_{\text{reg}}$ be the subset of regular elements in $\ka$,
i.e. those elements which belong to no hyperplane $\{\alpha=0\}$.
For $I$ a subset of $\kR^+$, let $$\ka^I:=\{x\in \ka \mid \forall
\alpha \in I,\ (\alpha,x)=0\}$$ be the face associated to $I$. Let
$\kR_I$ be the set of positive roots which are orthogonal to
$\ka^I$, and let $W_I$ be the subgroup of $W$ generated by the
$r_\alpha$ with $\alpha \in \kR_I$.
\newline
For $\xi \in \ka$, let $T_\xi$ be the Dunkl-Cherednik operator. It
is defined, for $f\in C^1(\ka)$, and $x\in \ka_{\text{reg}}$, by
$$T_\xi f(x)=\partial_\xi
f(x) + \sum_{\alpha \in \kR^+}k_\alpha
\frac{(\alpha,\xi)}{1-e^{-(\alpha,x)}}\{f(x)-f(r_\alpha x)
\}-(\rho,\xi)f(x),$$ where
$$\rho=\frac{1}{2}\sum_{\alpha \in \kR^+}k_\alpha \alpha.$$
The Dunkl-Cherednik operators form a commutative family of
differential-difference operators (see \cite{C} or \cite{O}). The
Heckman-Opdam Laplacian $\kL$ is defined by
$$\kL=\sum_{i=1}^{n} T_{\xi_i}^2,$$
where $\{\xi_1,\dots,\xi_n\}$ is any orthonormal basis of $\ka$
($\kL$ is independent of the chosen basis). Here is an explicit
expression (see the appendix), which holds for $f\in C^2(\ka)$ and
$x\in \ka_{\text{reg}}$:
\begin{eqnarray}
\label{explicit} \kL f(x) &=& \Delta f(x)+ \sum_{\alpha \in
\kR^+}k_\alpha \coth \frac{(\alpha,x)}{2}\partial_\alpha f(x)
+|\rho|^2f(x) \\
 &-& \nonumber \sum_{\alpha \in \kR^+}k_\alpha \frac{|\alpha|^2}{4\sinh^2 \frac{(\alpha,x)}{2}}
\{f(x)-f(r_\alpha x)\}
\end{eqnarray}
Let $\la \in \kh$. We denote by $F_\la$ the unique analytic
$W$-invariant function on $\ka$, which satisfies the differential
equations
$$p(T_\xi)F_\la=p(\la)F_\la \text{ for all W-invariant polynomials }p$$
and which is normalized by $F_\lambda(0)=1$ (in particular $\kL
F_\la=(\la,\la)F_\la$). We denote by $G_\lambda$ the unique
analytic function on $\ka$, which satisfies the differential and
difference equations
\begin{eqnarray}
\label{equations} T_\xi G_\la = (\la,\xi)G_\la \text{ for all }\xi
\in \ka,
\end{eqnarray}
and which is normalized by $G_\lambda(0)=1$.
\newline
\newline \textit{The $\textbf{c}$-function.} \newline We define the
function $\textbf{c}$ as follows (see \cite{H} or \cite{HO}):
$$\mathbf{c}(\la)=c_0 \prod_{\alpha \in
\kR^+}\frac{\Gamma(-(\la,\alpha^\vee)+\frac{1}{2}k_{\frac{\alpha}{2}})}
{\Gamma(-(\la,\alpha^\vee)+k_\alpha+\frac{1}{2}k_{\frac{\alpha}{2}})},$$
where $c_0$ is a positive constant chosen in such a way that
$\mathbf{c}(-\rho)=1$, and $k_{\frac{\alpha}{2}}=0$ if
$\frac{\alpha}{2} \notin \kR$. Observe that if
$$\mathbf{\pi} (\la):=\prod_{\alpha \in \kR_0^+}(\lambda,\alpha^\vee),$$ then the
function
$$\mathbf{b}(\la):=\mathbf{\pi}(\la)\mathbf{c}(\la),$$
is analytic in a neighborhood of $0$.
\begin{rem}\emph{
For the reader's convenience, let us point out a conventional
difference between our setting and symmetric spaces. There
$\Sigma$ denotes the root system and $m:\Sigma \to \N^*$ the
multiplicity function. Everything fits together if we set
$\kR=2\Sigma$ and $k_{2\alpha}=\frac{1}{2}m_\alpha$. Notice in
particular that $\rho$ is defined in the same way in both settings:
$$\rho =\frac{1}{2}\sum_{\alpha \in \Sigma^+} m_\alpha
\alpha=\frac{1}{2}\sum_{\alpha \in \kR^+}k_\alpha \alpha.$$ }
\end{rem}

\section{Estimates}

\label{secestim}
\subsection{Positivity and first estimates}
Let us begin with the following positivity result.
\begin{lem}
\label{positivite} Assume that $\la \in \ka$. Then the functions
$F_\lambda$ and $G_\lambda$ are real and strictly positive.
\end{lem}
\textbf{Proof of lemma:} Since
\begin{eqnarray}
\label{formFG} F_\lambda(x)=\frac{1}{|W|}\sum_{w\in W}G_\lambda
(w\cdot x), \text{   } x\in \ka,
\end{eqnarray}
it is enough to prove the lemma for $G_\la$. First of all, the
function $G_\la$ is real valued, since $G_\la$ and
$\overline{G_\la}$ satisfy the same equations (\ref{equations}),
and hence are equal. Assume next that $G_\la$ vanishes. Let $x$ be
a zero of $G_\la$ of minimal norm $r=|x|$. Consider first the case
where $x$ is a regular point, and take a vector $\xi$ in the same
chamber as $x$. As $G_\la$ is positive for $|x|<r$, we have
$$\partial_\xi G_\lambda(x) \le 0.$$
Writing down (\ref{equations}), we get
\begin{eqnarray}
\label{eqd}
\partial_\xi G_\la(x) &=& \sum_{\alpha \in \kR^+}k_\alpha
\frac{(\alpha,\xi)}{1-e^{-(\alpha,x)}}(G_\la(r_\alpha x)-G_\la(x))
+ (\rho+\la,\xi)G_\la(x).
\end{eqnarray}
Since for all roots $\alpha$,
$$\frac{(\alpha,\xi)}{1-e^{-(\alpha,x)}}\ge 0,$$
we deduce that $\partial_\xi G_\lambda(x)= 0$, and that
$G_\la(r_\alpha x) =0$ for every $\alpha \in \kR$. Hence $G_\la$
and $\nabla G_\la$ vanish at the point $x$ and furthermore at each
conjugate of $x$ under $W$. Differentiating (\ref{eqd}), we see
that every second order partial derivative of $G_\la$ vanishes on
the $W$-orbit of $x$. And similarly for all higher order
derivatives. Since $G_\la$ is analytic, we deduce that $G_\la
\equiv 0$. This contradicts the fact that $G_\la(0)=1$.
\newline Consider next the case where $x$ is singular and let $I=\{\alpha \in \kR^+ \mid (\alpha,x)=0 \}$. The equations (\ref{equations})
become now
\begin{eqnarray} \label{eqdiff} \partial_\xi
G_\la(x) &=& -\sum_{\alpha \in I}
2k_\alpha\frac{(\alpha,\xi)}{|\alpha|^2}\partial_\alpha G_\la(x) \\
                                &+& \nonumber \sum_{\alpha \in \kR^+\smallsetminus I} k_\alpha
\frac{(\alpha,\xi)}{1-e^{-(\alpha,x)}}(G_\la(r_\alpha
x)-G_\la(x))+(\rho+\la,\xi)G_\la(x).
\end{eqnarray}
We may argue as before, taking $\xi \in \ka^I$ in the same face as
$x$. Notice that the first sum vanishes in the right hand side of
(\ref{eqdiff}), and that
$$\partial_\xi(r_\alpha G_\la)(x)=\partial_{r_\alpha \xi}G_\la(r_\alpha x).$$
with $r_\alpha \xi$ in the same face as $r_\alpha x$. Eventually
we obtain that all partial derivatives of $G_\la$ along directions
belonging to $\ka^I$ vanish at $x$. Again since $G_\la$ is
analytic, it must vanish on $\ka^I$, which contradicts
$G_\la(0)=1$. This concludes the proof of the lemma. \hfill
$\blacksquare$
\newline
\newline
The next proposition is fundamental in order to have uniform
estimates in the parameter $\la \in \kh$.
\begin{prop}
\label{majorationbase} \begin{itemize} \item[(a)] For all $\la \in
\kh$,
$$|G_\la| \leq G_{\Re(\la)}.$$
\item[(b)] For all $\la \in \ka$ and for all $x\in \ka$
$$G_\la(x) \leq G_0(x)e^{\max_w(w\la,x)}.$$
\end{itemize}
\end{prop}
\textbf{Proof of the proposition:} For the first inequality, we
study the behavior of the ratio
$Q_\la=\frac{G_\la}{G_{\Re(\la)}}$. We must show that $|Q_\la|^2
\leq 1$. We will in fact prove that for all $\xi \in
\ka_{\text{reg}}$,
$$M(\xi,r):=\max_{w\in W} |Q_\la(rw \xi)|^2$$
is a decreasing function of $r\ge 0$. Since $M(\xi, 0)=1$ for all
$\xi$, the result will follow. First of all observe that the
function $M$ is continuous and right differentiable in the second
variable $r$. Then, using (\ref{equations}), we get
\begin{eqnarray*}
\partial_\xi |Q_\la|^2(x)=\sum_{\alpha \in \kR^+}\frac{2k_\alpha(\alpha,\xi)}{1-e^{-(\alpha,x)}}
(\Re\{Q_\la(x)\overline{Q_\la}(r_\alpha
x)\}-|Q_\la(x)|^2)\frac{G_{\Re(\la)}(r_\alpha
x)}{G_{\Re(\la)}(x)},
\end{eqnarray*}
for all $\xi$ and all $x$ regular. Hence if $x$ is a regular
element such that $$|Q_\la(x)|^2= \max_{w}|Q_\la(w x)|^2,$$ and if
$\xi$ is a positive multiple of $x$, we have
$$\partial_{\xi} |Q_\la|^2(x)\leq 0.$$
This means that
$$\frac{\partial M}{\partial r}(\xi,|x|) \leq 0,$$ where we consider right derivatives. So for every $\xi$ regular, and every
$r\ge 0$,
$$\frac{\partial M}{\partial r}(\xi,r) \leq 0.$$
In order to conclude, we need the following elementary lemma,
whose proof is left to the reader.
\begin{lem}
\label{lemana} Let $f:\R^+ \to \R$ be a continuous and right
derivable function. We denote by  $f'_d$ the right derivative of
$f$. If for all $x\in \R^+$, $f'_d(x)\le 0$, then $f$ is
decreasing.
\end{lem}
According to this lemma, we have $M(\xi,r)\le M(\xi,0)=1$, for all
$\xi\in \ka_{\text{reg}}$ and all $r\ge 0$. By continuity, this
inequality remains true if $\xi$ is singular. This concludes the
proof of the first inequality.
\newline
The second one is proved similarly, using the ratio
$$R_\la(x):=\frac{G_\la(x) e^{-\max_w(w\la,x)}}{G_0(x)}.$$
Specifically, if $x$ is regular and $\xi\in \ka$, then
\begin{eqnarray*}
\partial_\xi R_\la(x) &=& \sum_{\alpha
\in \kR^+}\frac{k_\alpha(\alpha,\xi)}{1-e^{-(\alpha,x)}}
(R_\la(r_\alpha x)-R_\la(x))\frac{G_0(r_\alpha x)}{G_0(x)} \\
  & +&
((\la,\xi)-\max_w(w\la,\xi))R_\la(x),
\end{eqnarray*}
where we consider again right derivatives. So if $x$ is such that
$$R_\la(x)=\max_w R_\la(w x)$$ and $\xi$ is a positive multiple
of $x$, then $$\partial_\xi R_\la(x)\le 0.$$ Therefore
$$N(\xi,r):=\max_{w\in W} R_\la(rw\cdot \xi)$$
is a decreasing function in $r\ge 0$, for all $\xi \in
\ka_{\text{reg}}$. We conclude as for the first inequality.
 \hfill
$\blacksquare$ \newline
\newline By averaging over the Weyl group, we deduce the following inequalities from Proposition
\ref{majorationbase}.
\begin{cor}\begin{enumerate}
\item
For all $\la \in \kh$,
$$|F_\la| \leq F_{\Re(\la)}.$$
\item For all $\la \in \ka$ and for all $x\in \ka$
$$F_\la(x) \leq F_0(x)e^{\max_{w\in W}(w\la,x)}.$$
\end{enumerate}
\end{cor}

\subsection{Local Harnack principles and sharp global estimates}

In this subsection we first establish two Harnack principles for
$G_\la$ and $F_\la$ when $\la \in \ka$, and next deduce sharp
global estimates of these functions $F_\la$ and of the function
$G_0$. Before stating the results we introduce some new notation.
Let $I$ be a subset of $\kR^+$, and let $d\le d'$ be two strictly
positive constants. We denote by $V^I(d,d')$ the following subset
of $\ka$:
$$V^I(d,d'):=\{x\in \ka \mid \forall \alpha \in \kR_I,\ |(\alpha,x)|\le d \text{ and }\forall \alpha
\notin \kR_I,\ |(\alpha,x)|> d' \}.$$ Let $x\in V^I(d,d')$, with
$I$ non empty. Let $p^I(x)$ denote its orthogonal projection on
$\ka^I$. Let $u\in \ka^I$ be such that for every $\alpha \notin
\kR_I$, $(\alpha,u)\text{sgn}((\alpha,x))\ge |\alpha|.$ Define now
the vectors $\xi_1(x)$, and $\eta_1(x)$ as follows:
\begin{eqnarray*}
\xi_1(x)=\frac{p^I(x)-x}{|p^I(x)-x|}+u, \text{ and  }
\eta_1(x)=\frac{p^I(x)-x}{|p^I(x)-x|}-u.
\end{eqnarray*}
We will sometime just write them $\xi_1$ and $\eta_1$ for simplify
the notation. Notice that everything was done in order that
\begin{eqnarray}
\label{xi} \forall \alpha \notin \kR_I,\
(\alpha,\xi_1(x))(\alpha,x)\ge0 \text{
 and  }(\alpha,\eta_1(x))(\alpha,x)\le0.
\end{eqnarray}
Naturally we have also
\begin{eqnarray}
\label{xideux} \forall \alpha \in \kR_I,\
(\alpha,\xi_1(x))(\alpha,x)=(\alpha,\eta_1(x))(\alpha,x)\le0.
\end{eqnarray}
 We denote by $p_1$ and
$q_1$ the projections of $x$ on $\ka^I$ along the directions
$\xi_1$ and $\eta_1$ respectively (we suppose that $d'$ is
sufficiently large in order that these projections still lie in
the same chamber than $x$). Then we denote by $p^\emptyset$ and
$q^\emptyset$ the orthogonal projections of $p_1$ and $q_1$
respectively on $V^\emptyset(d,d)$. We define also the vectors
$\xi_2$ and $\eta_2$ (like before we forget the dependence in $x$ in
the notation) by
\begin{eqnarray*}
\xi_2=\frac{p^\emptyset-p_1}{|p^\emptyset-p_1|}+u, \text{ and  }
\eta_2=\frac{q^\emptyset-q_1}{|q^\emptyset-q_1|}-u.
\end{eqnarray*}
Eventually let $p_2$ and $q_2$ be the projections on
$V^\emptyset(d,d)$ of $p_1$ and $q_1$ respectively along the
directions $\xi_2$ and $\eta_2$ (here again we suppose that $d'$
is sufficiently large in order that these projections lie in the
same chamber than $x$). We summarize these definitions in the
following figure \setlength{\unitlength}{.4cm}
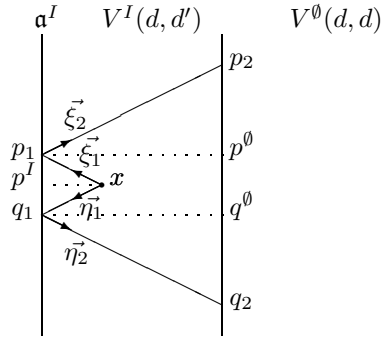
\begin{center}
\begin{picture}(12,12) \put(2,0){\line(0,1){10}}
\put(8,0){\line(0,1){10}} \put(4,5){\line(-2,1){2}}
\put(4,5){\line(-2,-1){2}} \put(2,6){\line(2,1){6}}
\put(2,4){\line(2,-1){6}} \put(4,5){\circle*{.2}}
\bezier{5}(4,5)(3,5)(2,5) \put(4,5){\vector(-2,1){1}}
\put(4,5){\vector(-2,-1){1}} \put(2,6){\vector(2,1){1}}
\put(2,4){\vector(2,-1){1}} \put(4.25,5){$x$} \put(1,4){$q_1$}
\put(1,6){$p_1$} \put(1,5){$p^I$} \put(4.25,5){$x$}
\put(8.25,6){$p^\emptyset$} \put(8.25,4){$q^\emptyset$}
\put(8.25,1){$q_2$} \put(8.25,9){$p_2$} \bezier{15}(2,6)(5,6)(8,6)
\bezier{15}(2,4)(5,4)(8,4) \put(1.75,10.25){$\ka^I$}
\put(4,10.25){$V^I(d,d')$} \put(10.25,10.25){$V^\emptyset(d,d)$}
\put(3.25,4){$\vec{\eta_1}$} \put(3.25,5.75){$\vec{\xi_1}$}
 \put(2.75,7){$\vec{\xi_2}$}  \put(2.75,2.5){$\vec{\eta_2}$}
\end{picture}
\end{center}
We can now state the lemma
\begin{lem}[Local Harnack principle 1]
\label{harnack1} Let $\la \in \ka$, and let $d$ and $d'$ be chosen
as above. There exist two constants $C>0$ and $c>0$ such that for
all $x \in V^I(d,d')$,
$$\max_{w \in W_I}G_\la(w x) \le C\min_{w \in
W_I}G_\la(wp_2(x)),$$ and
$$\min_{w \in W_I}G_\la(wx) \ge
c\max_{w \in W_I}G_\la(wq_2(x)). $$
\end{lem}
\textbf{Proof of the lemma:} We begin by the first inequality.
Let $x\in V^I$. First remark that $|x-p_1(x)|$ and $|x-q_1(x)|$
are bounded by a constant, say $h$, which depends only on $d$. We
introduce the function $M_\la$ defined on $\ka$ by:
$$M_\la(x)=\max_{w\in W_I} G_\la(w x).$$
Let $y$ be such that $G_\la(y)=M_\la(y).$ We have
\begin{eqnarray*}
\partial_{\xi_1} G_\la(y) &=& \sum_{\alpha \in \kR_I}k_\alpha
\frac{(\alpha,\xi_1)}{1-e^{-(\alpha,y)}}(G_\la(r_\alpha y)-G_\la(y))\\
                                       &+&\sum_{\alpha \in \kR^+\smallsetminus \kR_I}k_\alpha
\frac{(\alpha,\xi_1)}{1-e^{-(\alpha,y)}}(G_\la(r_\alpha
y)-G_\la(y))\\
                    & +&  (\rho+\la,\xi_1)G_\la(y)\\
                    &\ge& -\sum_{\alpha\in \kR^+\smallsetminus \kR_I}k_\alpha
\frac{(\alpha,\xi_1)}{1-e^{-(\alpha, y)}}G_\la(y)+ (\rho+\la,\xi_1)G_\la(y).\\
\end{eqnarray*}
The lower bound is deduced from our choice of $y$ and from the
properties of $\xi_1$ (\ref{xi}) and (\ref{xideux}). Now when
$\alpha \in \kR^+\smallsetminus \kR_I$, the ratio
$\frac{(\alpha,\xi_1)}{1-e^{-(\alpha,y)}}$ is bounded by a
constant which depends only on $d'$. Thus we can find a constant
$K$, which depends only on $d'$ and $\la$ such that for all $y\in
V^I(d,d')$,
$$\partial_{\xi_1} M_\la(y) \ge -K M_\la(y).$$ Here like in the proof of Proposition \ref{majorationbase}, we consider the right derivatives.
Still by Lemma \ref{lemana}, we get
\begin{eqnarray}
\label{mlambda} M_\la(x)\le e^{Kh}M_\la(p_1(x)). \end{eqnarray}
Now we introduce the function $N_\la$ defined on $\ka$ by
$$N_\la(x)=\min_{w \in W_I}G_\la(w x).$$
Observe already that $N_\la$ and $M_\la$ are equal on $\ka^I$, and
in particular in $p_1(x)$. Moreover, by the same technique as
above, we can find a strictly positive constant $K'$ such that
$$N_\la(p_1(x))\le e^{K'h}N_\la(p_2(x)).$$
Together with (\ref{mlambda}) this proves the first inequality of
the lemma. The second one can be proved exactly in the same way,
by using this time the intermediate point $q_1(x)$. \hfill
$\blacksquare$
\newline
\newline
We could deduce from this lemma a local Harnack principle for
$F_\la$ too. We will instead give a simple expression of the
gradient of $F_\la$, which implies such a principle. Moreover this expression will be needed in the proof of Theorem \ref{euler}.
\begin{lem}[Local Harnack principle 2]
\label{harnack} For all $x\in \overline{\ka_+}$ and for all $\la
\in \ka$, \begin{eqnarray} \label{nablaf} \nabla F_\la(x)=-
\frac{1}{|W|}\sum_{w\in W} w^{-1}(\rho-\la)\ G_\la(wx).
\end{eqnarray}
In particular,
$$|\nabla F_\la(x)| \le (|\rho|+|\la|)F_\la(x).$$
\end{lem}
\textbf{Proof of the lemma:} By differentiating (\ref{formFG}) we
get as above
$$\partial_\xi F_\la(x)=\frac{1}{|W|}\sum_{w\in W} \partial_{w\xi}G_\la(w x),$$
for all $\xi \in \ka$. Now we use the equations (\ref{equations}),
which gives
\begin{eqnarray*}
\partial_\xi F_\la(x) &=& \frac{1}{|W|}\sum_{w\in W}\sum_{\alpha \in \kR^+}k_\alpha \frac{(\alpha,w\xi)}{1-e^{-(\alpha,wx)}}\{G_\la(r_\alpha w x)-G_\la(w x)\}\\
                   &+& \frac{1}{|W|}\sum_{w\in W} (\rho+\la,w\xi)G_\la(w x)\\
                   &=& -\frac{1}{|W|}\sum_{w\in W}\sum_{\alpha \in \kR^+}k_\alpha (\alpha,w\xi)
                    \{\underbrace{\frac{1}{1-e^{-(\alpha,w\xi)}}+\frac{1}{1-e^{(\alpha,w\xi)}}}_{=1}\}G_\la(w x)\\
                   &+&  \frac{1}{|W|}\sum_{w\in W} (\rho+\la,w\xi)G_\la(w x)\\
                   &=& \frac{1}{|W|}\sum_{w\in W} (\la-\rho,w\xi)G_\la(w x).\\
\end{eqnarray*}
This proves the first claim of the lemma. The second one is an
easy consequence, using again (\ref{formFG}) and the positivity of
$G_\la$. \hfill $\blacksquare$
\newline
\newline
We can now deduce a sharp global estimate of $F_0$ which extends
the result of Anker \cite{A} to any multiplicities $k
> 0$. Recently Sawyer \cite{Saw} has obtained the same result for root
systems of type $A$, using explicit formulas.
\begin{theo}
\label{estimate1} In $\overline{\ka_+}$,
$$F_0(x) \asymp e^{-(\rho,x)}\prod_{\alpha \in \kR_0^+} (1+(\alpha,x)).$$
\end{theo}
\textbf{Proof of the theorem:} We resume the proof in \cite{A},
that we sketch. The local Harnack principle for $F_0$ (which was
deduced in \cite{A} from Harish-Chandra's integral formula) allows
us to move the estimate away from the walls in $\ka_+$. There we
expand $F_0$, using the Harish-Chandra series
\begin{eqnarray*} F_\la(x)= \sum_{w\in W} \sum_{q \in Q^+}
\mathbf{c}(w\la) \Gamma_q(w\la)e^{(w\la-\rho-q,x)} \end{eqnarray*}
that we multiply by $\mathbf{\pi}(\la)$ in order to remove the
singularity of the $\mathbf{c}$-function at the origin. Then we
differentiate with respect to
$\mathbf{\pi}(\frac{\partial}{\partial \la})|_{\la=0}$, in order
to recover $F_0$, up to a positive constant. As a result we obtain
a converging series
$$F_0(x)=\sum_{q\in Q^+} F_{q}(x)e^{-(\rho+q,x)},$$ with
polynomial coefficients $F_{q}$ and leading term
$$F_{0}e^{-(\rho,x)}\sim \text{const.}
\mathbf{\pi}(x)e^{-(\rho,x)}.$$
 \hfill $\blacksquare$
\begin{rem}\emph{
We may estimate in a similar way the function $F_\la$ when $\la$
is real. The result reads as follows: for any $\la \in
\overline{\ka_+}$,
$$F_\la(x) \asymp \prod_{\alpha \in \kR_0^+ \mid (\alpha,\la)=0}
(1+(\alpha,x))e^{(\la-\rho,x)}$$ on $\overline{\ka_+}$.}
\end{rem}
Let us turn to the function $G_0$. For $x\in \ka$, we denote by
$x^+$ its unique symmetric in $\overline{\ka_+}$.
\begin{theo}
\label{estimate2} In $\ka$, \begin{eqnarray} \label{estimationG0}
G_0(x) \asymp \prod_{\alpha \in \kR_0^+ \mid (\alpha,x)\ge
0}(1+(\alpha,x))e^{(-\rho,x^+)}.
\end{eqnarray}
\end{theo}
\textbf{Proof of the theorem:} Let us first show that $G_\la$ has
a series expansion in each chamber, like it was done by Opdam in
the negative chamber $\ka_-$ \cite{O}. We resume his proof. He
first obtained that there exists a polynomial $p$ such that for
all $x\in \ka_{\text{reg}}$,
$$(\prod_{\alpha \in \kR_0^+}(\la,\alpha^\vee)-k_\alpha-2k_{2\alpha})G_\la(x)=p(\la,T_\xi)F_\la(x).$$
By expanding $F_\la$ and $T=T_\xi$ in each chamber, we
find developments of the function $G_\la$:
$$G_\la(x)= \sum_{w'\in W}\mathbf{c}(w^{-1}w'\la)\sum_{q\in w Q^+}
G^{w,w'}_{\la,q}e^{(w'\la-w\rho-q,x)}$$ for all $x\in w\ka_+$.
Moreover Opdam has proved that $G^{w_0,w'}_{\la,0}$ is equal to
$|W|\delta_{1,w'}\mathbf{\mathbf{\pi}(\la)},$ where $w_0$ denotes
the longest element in $W$. Now we apply the same technique as in
Theorem \ref{estimate1}. First we multiply these developments by
$\mathbf{\pi}(\la)$, and then we differentiate with respect to
$\mathbf{\pi}(\frac{\partial}{\partial \la})|_{\la=0}$. We get
developments of the function $G_0$ in each chamber:
\begin{eqnarray}
\label{devserieG0}
G_0(x)=\sum_{q\in w Q^+} G^w_{q}(x)e^{-(w\rho+q,x)}
\end{eqnarray}
for all $x\in w\ka_+$, where the $G^w_{q}$ are real polynomials.
Moreover according to the above mentioned result of Opdam, we see that
$G^{w_0}_{0}$ is a strictly positive constant. Recall some basic notation. The length $l(w)$ of an
element of $W$ is defined by
$$l(w)=|\kR_0^+\cap w\kR_0^-|.$$
Recall that $\Pi$ denotes the set of simple roots in $\kR^+$. Each
$q\in Q^+$ writes $q=\sum_{\alpha \in \Pi} n_\alpha \alpha$, with
$n_\alpha \in \N$. We denote by $|q|:=\sum_{\alpha \in \Pi}
n_\alpha$ the length of $q$. For $q'\in Q^+$, we write $q'\le q$,
if $q-q'\in Q^+$. Naturally we have similar definitions on $wQ^+$,
where we denote by $|q|_w$ the length of any $q\in wQ^+$ and we
write $q' \le_w q$, if $q' \in wQ^+$ and $q-q'\in wQ^+$. Consider
the polynomials
$$\pi_w(x)=\prod_{\alpha \in \kR^+_0\cap
w\kR_0^+}(\alpha^\vee,x)\quad \text{and }\quad
\tilde{\pi}_w(x)=\prod_{\alpha \in \kR^+_0\cap
w\kR_0^+}\Big(1+(\alpha^\vee,x)\Big).$$ We need the following
lemma, which will be used throughout the proof of Theorem \ref{estimate2}.
\begin{lem}
\label{cornerlemma} Let $w\in W$. \begin{enumerate}
\item If $\alpha \in \Pi \cap w\kR^+$, then $\pi_{r_\alpha w}(r_\alpha x)=
\frac{\pi_w(x)}{(\alpha^\vee,x)}$, for all $x\in
\ka_{\text{reg}}$.
\item If $\alpha \in \kR_0^+\cap w\kR_0^+$, then
$\tilde{\pi}_{r_\alpha w}(r_\alpha x)\le
\frac{\tilde{\pi}_w(x)}{1+(\alpha^\vee,x)}$, for all $x\in
w\ka_+$.
\item If $\alpha \in \kR_0^-\cap w\kR_0^+$, then there exists a constant $C>0$, such that
$\tilde{\pi}_{r_\alpha w}(r_\alpha x)\le C
\tilde{\pi}_w(x)(1+(\alpha^\vee,x))^{|\kR^+|}$, for all $x\in
w\ka_+$.
\end{enumerate}
\end{lem}
\textbf{Proof of the lemma:} Let us prove the first claim. Since
$\alpha \in \Pi$, $r_\alpha$ maps $\kR_0^+\smallsetminus
\{\alpha\}$ onto itself, hence $\kR_0^+\cap r_\alpha w\kR_0^+$
onto $(\kR_0^+\cap w\kR_0^+)\smallsetminus \{\alpha\}$. The first
claim follows. \newline Let us prove the second claim. We define
therefore an injective map $i$ from $\kR_0^+\cap r_\alpha
w\kR_0^+$ into $(\kR_0^+\cap w\kR_0^+)\smallsetminus \{\alpha\}$,
such that $r_\alpha \beta \le_w i(\beta)$ for all $\beta$. The
second claim will follow. Let $\beta \in \kR_0^+\cap r_\alpha
w\kR_0^+$. If $r_\alpha \beta \in \kR_0^+$, then we set
$i(\beta)=r_\alpha \beta$. Otherwise, we have $(\alpha,\beta)\ge
0$. Hence $r_\alpha \beta \le_w \beta$. But $r_\alpha \beta \in
w\kR_0^+$, and therefore $r_\alpha \beta \ge_w 0$. Thus $\beta \in
\kR_0^+\cap w\kR_0^+$ and we set $i(\beta)=\beta$. The map $i$
defined this way has all required properties.
\newline Let us prove the third claim. We define this time an injective map $i$
from $I\subset \kR_0^+\cap r_\alpha w\kR_0^+$ into $\kR_0^+\cap
w\kR_0^+$ such that, if $\beta \in I$, then $r_\alpha \beta\le_w
i(\beta)+|(\alpha^\vee,\beta)|\alpha$, and otherwise $r_\alpha
\beta \le_w |(\alpha^\vee,\beta)|\alpha$. The third claim will
follow. Assume that $\beta \in \kR_0^+\cap r_\alpha w\kR_0^+$. If
$r_\alpha \beta \in \kR_0^+$, then we set $i(\beta)=r_\alpha
\beta$. Otherwise $(\alpha,\beta) \le 0$. Next, either $\beta
\in w\kR_0^+$, in which case $r_\alpha \beta \le_w
\beta+|(\alpha^\vee,\beta)|\alpha$, and we set
$i(\beta)=\beta$. Or $\beta \in w\kR_0^-$ in which case $r_\alpha \beta\le_w
|(\alpha^\vee,\beta)|\alpha$. The map $i$ defined this way has
all required properties. \hfill $\blacksquare$
\newline
\newline
By expanding $G_0$ in (\ref{eqd}) according to (\ref{devserieG0}) we get
\begin{eqnarray}
\label{formulerecurrence} \nabla G^w_{q}(x)&=& G^w_{q}(x)q+\sum_{\alpha \in
\kR^+\cap w\kR^+}k_\alpha G^{r_\alpha
w}_{r_\alpha q}(r_\alpha x)\alpha\\
&+&\nonumber \sum_{\alpha \in w\kR^+}k_\alpha\sum_{j\in \N^*}\{G^{r_\alpha
w}_{r_\alpha (q-j\alpha)}(r_\alpha x)-G^w_{(q-j\alpha)}(x)\}\alpha,
\end{eqnarray}
for all $w\in W$, all $q\in wQ^+$, and all $x\in w\overline{\ka_+}$.
\newline
\newline
\textit{Step $1$:} Let us first establish the estimate
$$|G^w_{0}(x)|\le
C\tilde{\pi}_w(x)\quad \forall w\in W,\quad \forall x\in w\overline{\ka_+}.$$
It is obvious for $w=w_0$. Let us prove it by induction on $l(w_0)-l(w)$.
For $q=0$, (\ref{formulerecurrence}) amounts to
$$\partial_\xi G^w_{0}(x)=\sum_{\alpha \in
\kR^+\cap w\kR^+}k_\alpha (\alpha,\xi) G^{r_\alpha
w}_{0}(r_\alpha x).$$
Using the induction hypothesis and Lemma \ref{cornerlemma}, we get
\begin{eqnarray*}
\partial_\xi G^w_{0}(x)&\le & C \sum_{\alpha \in
\kR^+\cap w\kR^+}k_\alpha (\alpha,\xi) \tilde{\pi}_{r_\alpha
w}(r_\alpha x)\\
&\le & C\sum_{\alpha \in
\kR^+\cap w\kR^+}k_\alpha \frac{(\alpha,\xi)}{1+(\alpha^\vee,x)}\tilde{\pi}_w(x)\\
&=&C\partial_\xi \tilde{\pi}_w(x)
\end{eqnarray*}
for all $x\in w\overline{\ka_+}$ and $\xi\in w\overline{\ka_+}$,
in particular for $\xi \in \R^+x$. Since $G^w_{0}(0)\le C$
provided $C$ is large enough, we obtain the upper estimate $$G^w_{0}(x)\le
C\tilde{\pi}_w(x)\quad \forall x\in w\overline{\ka_+}.$$ The same argument
yields the lower estimate
$$G^w_{0}(x)\ge
-C\tilde{\pi}_w(x)\quad \forall x\in w\overline{\ka_+}.$$
\textit{Step $2$:} Let us next establish the following estimate:
There exist a constant $C>0$ and $h\in \ka_+$, such that for every
$w\in W$, $q\in wQ^+$ and $x\in C^w_h:=wh+w\overline{\ka_+}$,
\begin{eqnarray}
\label{recurrenceq} |G^w_{q}(x)|\le
C^{|q|_w}\tilde{\pi}_w(x)(1+q(x))^{|\kR^+|}.
 \end{eqnarray}
The case $q=0$ was considered in step $1$. Let $q\in
Q^+\smallsetminus \{0\}$ and $w\in W$. Assume that
(\ref{recurrenceq}) holds for all $(q',w')\in w'Q^+\times W$ such
that $|q'|_w< |q|_w$ or such that $|q'|_w=|q|_w$ and $l(w')<l(w)$.
Using (\ref{formulerecurrence}), the induction hypothesis and
Lemma \ref{cornerlemma}, we get
\begin{eqnarray} \label{majorationq}
\partial_\xi\Big[C^{|q|_w}\tilde{\pi}_w
(1+q)^{|\kR^+|}-G^w_{q}\Big](x)\ge (q,\xi)\Big[|\kR^+|
C^{|q|_w}\tilde{\pi}_w (1+q)^{|\kR^+|-1}- G_{q}^w\Big](x),
\end{eqnarray} for all $\xi \in w\ka_+$ and all $x\in wh+w\ka_+$,
provided $C>0$ is large enough. Using now (\ref{devserieG0}) at
the point $wh$ we can also assume, by taking again larger $C$ if
necessary, that $$G_{q}^w(wh)\le C^{|q|_w},$$ for all $q\in wQ^+$.
Let now $u \in wh+w\overline{\ka_+}$ be such that
$(1+(q,u))=|\kR^+|$. Equation (\ref{majorationq}) implies that
\begin{eqnarray}
\label{majorationqbis}[C^{|q|_w}\tilde{\pi}_w
(1+q)^{|\kR^+|}-G^w_{q}](x)\ge 0,
\end{eqnarray} for all $x$ in the segment $[wh,u]$. For $x=wh+t(u-wh)$ with $t\ge
1$, we have
$$\partial_{u}\Big[C^{|q|_w}\tilde{\pi}_w
(1+q)^{|\kR^+}-G^w_{q}\Big](x)\ge
(q,u)\frac{|\kR^+|}{(1+(q,x))}\Big[ C^{|q|_w}\tilde{\pi}_w
(1+q)^{|\kR^+|}- G_{q}^w\vee 0\Big](x).$$ Thus
(\ref{majorationqbis}) holds also for $x=wh+t(u-wh)$ with $t\ge1$.
This proves the upper estimate $$G^w_{q}(x)\le
C^{|q|_w}\tilde{\pi}_w(x)(1+q(x))^{|\kR^+|},$$ in $C^w_h$. The
same argument gives the lower estimate
$$G^w_{q}(x)\ge
-C^{|q|_w}\tilde{\pi}_w(x)(1+q(x))^{|\kR^+|}.$$ \textit{Step $3$:}
Let us now find a lower bound for $G^w_{0}$. We prove by induction
on $l(w_0)-l(w)$ that there exist a constant $c>0$ and $h\in
\ka_+$, such that $$G^w_{0}(x)\ge c \pi_w(x)$$ for all $x \in
C_h^w$. We suppose that it is true for $w'$ such that $l(w')>l$
and we consider $w$ of length $l$. By the induction hypothesis
there exists some $h\in \ka_+$ and $c>0$ such that, $G^{r_\alpha
w}_0(r_\alpha x)\ge c\pi_{r_\alpha w}(r_\alpha x)$, for all $x\in
C_h^w$ and all $\alpha \in \kR^+\cap w\kR^+$. Let now $c'>0$ be
another constant. Assume that for some $x_0 \in C_h^w$,
$$[G^w_{0}-c'\pi_w](x_0)\le [G^w_{0}-c'\pi_w](wh)-1,$$ and suppose that $x_0$
is such element of minimal norm in $C_h^w$. Let
$(\alpha^*)_{\alpha \in w\Pi}$ be the dual basis of $w\Pi$, i.e.
for $\alpha$ and $\beta$ in $w\Pi$, $(\alpha^*,\beta)=0$ if
$\alpha \neq \beta$ and $=1$ otherwise. Let $\alpha_0 \in w\Pi$ be
such that $(\alpha_0,x_0-h)=\max_{\beta \in w\Pi} (\beta,x_0-h)$.
It implies that, for small $\epsilon>0$ at least, $x_0-\epsilon
\alpha_0^* \in C_h^w$. Hence
$$\partial_{\alpha_0^*}[G^w_{0}-c'\pi_w](x_0)\le 0.$$ On the
other hand we know that for $x\in w\ka_+$,
\begin{eqnarray}
\label{formulederivees} \nabla [G^w_{0}-c'\pi_w](x)= \sum_{\beta
\in \kR^+\cap w\kR^+} \beta[k_\beta G^{r_\beta w}_0(r_\beta
x)-\frac{2c'}{|\beta|^2}\frac{\pi_w(x)}{(\beta^\vee,x)}].
\end{eqnarray} Now we need the following elementary lemma.
\begin{lem}
\label{lemsimple} Let $\alpha \in w\Pi$. Assume that there exists
$\beta \in \kR_0^+\cap w\kR^+$, such that $\alpha \le_w \beta$.
Then there exists $\gamma \in \Pi \cap w\kR^+$, such that $\alpha
\le_w \gamma$.
\end{lem}
\textbf{Proof of the lemma:} Let $\beta=\sum_{\gamma \in \Pi}
n_\gamma \gamma$ be the decomposition of $\beta$ in $\Pi$. Since
$\beta \in w\kR_0^+$, there exists $\gamma\in \Pi\cap w\kR_0^+$
such that $n_\gamma>0$. We see moreover that
$\tilde{\gamma}:=\sum_{\gamma \in \Pi\cap w\kR_0^+} n_\gamma
\gamma \in w\kR^+$, and that $\beta \le_w \tilde{\gamma}$, which
concludes the proof of the lemma. \hfill $\blacksquare$
\newline
\newline
Suppose now that there does not exist $\gamma\in \Pi\cap w\kR^+$
such that $\alpha_0\le_w \gamma$. Then by Lemma \ref{lemsimple},
no other $\beta \in \kR^+\cap w\kR^+$ satisfies $\alpha_0\le_w
\beta$. Thus from equation (\ref{formulederivees}) we get that for
all $y$ in the segment between $x_0-\epsilon \alpha_0^*$ and
$x_0$, $\partial_{\alpha_0^*}[G^w_{0}-c'\pi_w](y)=0$, which
contradicts the initial hypothesis on $x_0$. We conclude that
there exists $\gamma \in \Pi \cap w\kR^+$ such that $\alpha_0\le_w
\gamma$. Again from (\ref{formulederivees}) we get
$$\partial_{\alpha_0^*}[G^w_{0}-c'\pi_w](x_0)\ge c_1c\pi_{r_\gamma
w}(r_\gamma x_0)-c_2c'\frac{\pi_w(x_0)}{(\alpha_0^\vee,x_0)},$$
where $c_1$ and $c_2$ are positive constants. But with the first
point of Lemma \ref{cornerlemma} we have $\pi_{r_\gamma
w}(r_\gamma x_0)=\frac{\pi_w(x_0)}{(\gamma^\vee,x_0)}$. Moreover
by our choice of $\alpha_0$, we have $(\gamma,x_0)\le
|\gamma|_w(\alpha_0,x_0)$. Thus if $c'$ is sufficiently small we
get $$\partial_{\alpha^*_0}[G^w_{0}-c'\pi_w](x_0)>0$$ and a
contradiction. The induction hypothesis for $w$ follows.
\newline
\newline
Putting now the third steps together, we get the desired
estimate of $G_0$ away from the walls. With Lemma \ref{harnack1},
this concludes the proof of the theorem. \hfill $\blacksquare$
\newline
\newline
The preceding theorem has for us a very important consequence. Let
$E$ be the Euler operator. It is defined for $f$ regular, and
$x\in \ka$, by $Ef(x)=(x,\nabla f(x))$. The following theorem
generalizes the analogue result of \cite{ABJ} in the setting of
symmetric spaces. Our proof is in a certain sense more elementary
than in \cite{ABJ}, because we do not make use of the descent
technique of Harish-Chandra. \newline The first claim of the
theorem will be needed in the estimate of the heat semigroup
(Proposition \ref{estimchaleur}). It will also be used in the
study of the asymptotic convergence of the $F_0$-processes (see
\cite{ABJ}). It will allow us in \cite{Schap} to generalize some
results of Anker, Bougerol, and Jeulin \cite{ABJ} for all $k>0$.
The second claim is just a technical result needed in the proof of
the estimate of the heat kernel (see Theorem \ref{estimchaleur}).
\begin{theo}
\label{euler} \begin{enumerate} \item There exists a constant
$K>0$ such that for any $x \in \overline{\ka_+}$,
$$0\le E[\log (e^\rho F_0)](x)\le K.$$
\item We have the two following estimates
\begin{eqnarray*}
E[\log (e^\rho F_0)](x)=|\kR_0^+| + \kO(\frac{1}{1+\min_{\alpha
\in
\kR^+}(\alpha,x)}),\\
\sum_{\alpha \in
\kR_0^+}\frac{(\alpha,x)}{\sqrt{1+(\alpha,x)^2}}\partial_\alpha
(\log (e^\rho F_0))(x) \asymp \frac{1}{1+\min_{\alpha \in
\kR^+}(\alpha,x)}.
\end{eqnarray*}
\end{enumerate}
\end{theo}
\textbf{Proof of the theorem:} With the formula (\ref{formFG})
and (\ref{nablaf}), we get for any $x \in \overline{\ka_+}$,
\begin{eqnarray} \label{formuleeuler} E[\log(e^\rho
F_0)](x)=\frac{1}{|W|} \sum_{w\in W}[(\rho,x)-(\rho,w
x)]\frac{G_0(w x)}{F_0(x)}.
\end{eqnarray}
\begin{enumerate}
\item Formula (\ref{formuleeuler}) proves already the first inequality.
For the second inequality we show by induction on the length
$l(w)$ of $w\in W$ that for all $x \in \overline{\ka_+}$,
\begin{eqnarray}
\label{recurrence} (\rho,x)-(\rho,w x)\le K'l(w)\max_{\alpha \in
\kR_0^+\cap w\kR^-}|(\alpha,w x)|,
\end{eqnarray}
where $K'=\max_{\alpha \in \kR_0^+} (\rho,\alpha^\vee)$ is a
constant. Suppose that the induction hypothesis is true for all
$w$ of length less or equal to $l$. Let $v\in W$ be of length
$l+1$. Let $\alpha \in \Pi \cap v\kR^-$, and let $w=r_\alpha v$.
We have $l(w)=l$. Moreover since $\alpha \in \Pi$, $r_\alpha$ maps
$\kR_0^+\cap w\kR_0^-$ onto $(\kR_0^+\cap
v\kR_0^-)\smallsetminus\{\alpha \}$. But for all $x\in \ka$,
$$(\rho,x)-(\rho,vx)= (\rho,x)-(\rho,wx)-(\alpha,v x)(\rho,\alpha^\vee).$$
Thus (\ref{recurrence}) follows for $v$ by using the induction
hypothesis. Now with
Theorems \ref{estimate1} and \ref{estimate2}, the first claim
is proved.
\item These estimates result also from Formula
(\ref{formuleeuler}) and the global estimates (Theorems
\ref{estimate1} and \ref{estimate2}) of $G_0$ and $F_0$. The fact
that $|\kR_0^+|$ is the limit of $E[\log(e^\rho F_0)](x)$ when
$(\alpha,x)\to \infty$ for all $\alpha$ can be seen exactly like
in \cite{ABJ} by expanding the functions $F_\la$ in series. This
finishes the proof of the theorem. \hfill $\blacksquare$
\end{enumerate}

\subsection{Estimates of the derivatives}

In this subsection we estimate the derivatives of the
hypergeometric function $G_\la(x)$, first in $x$ alone and next
jointly in $(\la,x)$.
\begin{prop}
\label{derivx} Let $p$ be a polynomial of degree $N$. Then there
exists a constant $C$ such that, for any $\la \in \kh$ and for any
$x\in\ka$,
$$|p(\frac{\partial}{\partial x})G_\la(x)| \leq C(1+|\la|)^N F_0(x)e^{\max_w \Re
(w\la,x)}.$$
\end{prop}
\textbf{Proof of the proposition:} According to Proposition
\ref{majorationbase}, we know that this estimate holds with no
derivative.
\newline
\textit{Step $1$ : Estimate away from walls}
\newline
By induction, Formula (\ref{eqd}) allows us to express on
$\ka_{\text{reg}}$ derivatives of $G_\la$ in terms of lower order
derivatives and to estimate them away from walls. More precisely
we obtain this way the desired estimate when $x$ stays at distance
$\ge \frac{\epsilon}{1+|\la|}$ from walls.
\newline
\textit{Step $2$ : Estimate on faces} \newline Assume that $x$
lies in a face $\ka^I$ (of minimal dimension), then (\ref{eqd})
becomes (\ref{eqdiff}), which writes also
\begin{eqnarray}
\label{eqdege}
\partial_{A_{w,I}(\xi)}G_\la(x)&=&\sum_{\alpha \in \kR^+\smallsetminus
\kR_I}k_\alpha \frac{(\alpha,\xi)}{1-e^{-(\alpha,x)}}(G_\la
(r_\alpha x)-G_\la(x))\\
&+&\nonumber (\rho+\la,\xi)G_\la(x),
\end{eqnarray}
where $$A_{w,I}(\xi)=\xi + 2 \sum_{\alpha \in
\kR_I}\frac{k_\alpha}{|\alpha|^2}(\alpha,\xi)\alpha.$$ Notice that
the linear map $A_{w,I}:\ka\to\ka$ is one-to-one, since the
expression
$$(A_{w,I}(\xi),\xi)=|\xi|^2+2\sum_{\alpha\in
\kR_I}\frac{k_\alpha}{|\alpha|^2}(\alpha,\xi)^2$$ is strictly
positive for all nonzero $\xi$. By induction, (\ref{eqdege})
yields the following estimate: for every $\epsilon>0$, there
exists a constant $C\ge 0$ such that, for all multi-indices
$\kappa$, for all $\la \in \kh$ and for $x\in \ka^I$ such that
$\min_{\alpha \in \kR^+ \smallsetminus \kR_I} |(\alpha,x)|\ge
\frac{\epsilon}{1+|\la|}$,
\begin{eqnarray}
\label{derivface} |(\frac{\partial}{\partial x})^\kappa
G_\la(x)|\le |\kappa|! C^{|\kappa|}(1+|\la|)^{|\kappa|}
F_0(x)e^{\max_{w\in W}(w\Re \la,x)}.
\end{eqnarray}
\newline
\textit{Step $3$ : Estimate near the faces} \newline If $x$ is
near a face $\ka^I$, we use (\ref{derivface}) and the Taylor
development of $G_\la$ in the orthogonal projection of $x$ on
$\ka^I$. More precisely let $\epsilon>0$ be such that
$C\epsilon<1$, where $C$ is the constant appearing in
(\ref{derivface}). Then there exists a constant $C'>0$ such that,
for all multi-indices $\kappa$, for all $\la\in \kh$ and for $x\in
\ka$ at distance $\le \frac{\epsilon}{1+|\la|}$ from $\ka^I$, such
that $\min_{\alpha \in \kR^+ \smallsetminus \kR_I} |(\alpha,x)|\ge
\frac{\epsilon}{1+|\la|}$,
\begin{eqnarray}
\label{derivnear}|(\frac{\partial}{\partial x})^\kappa
G_\la(x)|\le C'(1+|\la|)^{|\kappa|} F_0(x)e^{\max_{w\in W}(w\Re
\la,x)}.
\end{eqnarray}
\newline
\textit{Step $4$ : Conclusion} \newline Now we first use the step
$3$ near the origin. We get $\epsilon_0>0$ and $C_0>0$, such that
(\ref{derivnear}) holds (with $C_0$ in place of $C'$) for $x\in
\ka$ at distance $\le \frac{\epsilon_0}{1+|\la|}$ from the origin.
Then we use the step $3$ near the faces of dimension $1$. We get
$\epsilon_1$ and $C_1$ such that (\ref{derivnear}) holds for $x\in
\ka$ at distance $\le \frac{\epsilon_1}{1+|\la|}$ from any face of
dimension $1$, and at distance $\ge \frac{\epsilon_0}{1+|\la|}$
from the origin. And like this we get successively, for each $d\in
\N$, constants $\epsilon_d>0$ and $C_d$ associated to the faces of
dimension $d$. Eventually we conclude with the first step. \hfill
$\blacksquare$
\newline
\newline
We can now derive the fundamental estimate:
\begin{theo}
\label{estimderivees} Let $p$ and $q$ be polynomials of degree $M$
and $N$. Then there exists a constant $C$ such that, for all $\la
\in \kh$ and for all $x\in\ka$,
$$|p(\frac{\partial}{\partial \la})q(\frac{\partial}{\partial x})G_\la(x)| \leq C(1+|x|)^M(1+|\la|)^NF_0(x)e^{\max_w \Re
(w\la,x)}.$$
\end{theo}
\textbf{Proof of the Theorem:} The proof is standard. Theorem
\ref{estimderivees} is deduced from Proposition \ref{derivx} using
Cauchy's formula. More precisely, one integrates $G_\la(x)$ in the
variable $\la$ over $n$-tori with radii comparable to
$\frac{1}{1+|x|}$. \hfill $\blacksquare$
\begin{rem} \emph{This estimate holds true for $F_\la$ too.}
\end{rem}

\section{Hypergeometric Fourier transform and Schwartz spaces}

We first recall the definitions of the hypergeometric Fourier
transform and of its inverse, according to Cherednik \cite{C2}.
Let $\mu$ be the measure on $\ka$ given by
$$d\mu(x)=\underbrace{\prod_{\alpha \in \kR^+}|2\sinh(\frac{\alpha}{2},x)|^{2k_\alpha}}_{:=\delta(x)}dx.$$ The
hypergeometric Fourier transform $\kH$ is defined for nice
functions $f$ on $\ka$ by
\begin{eqnarray}
\label{harish}
 \kH(f)(\la)=\int_\ka f(x)G_{\la}(-x)d\mu(x), \text{  }\forall \la \in \kh. \end{eqnarray} Let $\nu$ be the
asymmetric Plancherel measure on $i\ka$ defined by
$$d\nu
(\la)=c\prod_{\alpha \in \kR^+}
\frac{\Gamma((\la,\alpha^\vee)+k_\alpha+\frac{1}{2}k_{\frac{\alpha}{2}})\Gamma(-(\la,\alpha^\vee)+k_\alpha+\frac{1}{2}k_{\frac{\alpha}{2}}+1)}{\Gamma((\la,\alpha^\vee)+\frac{1}{2}k_{\frac{\alpha}{2}})\Gamma(-(\la,\alpha^\vee)+\frac{1}{2}k_{\frac{\alpha}{2}}+1)}d\la,$$
where $c$ is a normalizing constant. The inverse transform $\kI$
is given for nice functions $h$ by
\begin{eqnarray}
\label{forminv}
\kI(h)(x)=\int_{i\ka}h(\la)G_\la(x)d\nu(\la),\text{ }\forall x\in
\ka.
\end{eqnarray}
In the case $k=0$, $\kH$ and $\kI$ reduce to the classical
Euclidean Fourier transform
$$\kF(f)(\la)=\int_\ka f(x)e^{-(\la,x)}dx$$
and its inverse
$$\kF^{-1}(h)(x)=(2\pi)^{-n}\int_{i\ka}h(\la)e^{(\la,x)}d\la.$$
We shall consider the following function spaces. The classical
Schwartz space on $i\ka$ is denoted by $\kS (i\ka)$. Its topology
is defined by the semi-norms
$$\tau_{p,N}(h)=\sup_{\la \in
i\ka}(1+|\la|)^N|p(\frac{\partial}{\partial \la})h(\la)|,$$ where
$p$ is any polynomial and $N\in \N$. As usual $C_c^\infty(\ka)$
denotes the space of $C^\infty$ functions on $\ka$ with compact
support and $C_{\Gamma}^\infty(\ka)$ the subspace of functions
with support in a given compact subset $\Gamma$. Let us denote by
$\kC(\ka)$ the space of $C^\infty$ functions on $\ka$, such that
for all polynomials $p$ and all $N\in \N$,
$$\sup_{x\in \ka} (1+|x|)^NF_0(x)^{-1}|p(\frac{\partial}{\partial x}) f(x)| <
+\infty,$$ It is the Schwartz space on $\ka$ associated to the
measure $\mu$. Its topology is defined by the semi-norms
$$\sigma_{p,N}(f)=\sup_{x\in
\ka}(1+|x|)^NF_0(x)^{-1}|p(\frac{\partial}{\partial x})f(x)|.$$
Notice that according to Proposition \ref{estimate1}, we may
replace $F_0(x)$ by $e^{-(\rho,x^+)}$ in the definition of
$\kC(\ka)$ and its topology. Let us recall that $x^+$ is the only
point in the orbit $W\cdot x$ which lies in $\overline{\ka_+}$.
\begin{lem}
\begin{enumerate}
\item $\kC(\ka)$ is a Fréchet space.
\item $C_c^\infty(\ka)$ is a dense subspace of $\kC(\ka)$.
\end{enumerate}
\end{lem}
\textbf{Proof of the lemma:} These facts are standard. The second
one is proved for example in \cite{D}, more precisely in Appendix
$A$ by M. Tinfou. \hfill $\blacksquare$
\newline
\newline
Eventually, the Paley-Wiener space $P W(\kh)$ consists of all
entire functions $h$ on $\kh$ which satisfy the following growth
condition:
$$ \exists R\ge 0,\ \forall N\in \N,\ \sup_{\la \in \kh}
(1+|\la|)^Ne^{-R|\Re \la |} h(\la) <\infty.$$ Given a $W$-invariant
convex compact subset $\Gamma$ in $\ka$, $P W_{\Gamma}(\kh)$
denotes the subspace of $P W(\kh)$ defined by the specific
condition
$$\forall N\in \N,\ \sup_{\la \in \kh}
(1+|\la|)^Ne^{-\gamma(-\Re \la)} h(\la) <\infty.$$ Here
$\gamma(\la)=\sup_{x\in \Gamma} (\la,x)$ is the gauge associated
to the polar of $\Gamma$.
\newline
The mapping properties of the hypergeometric Fourier transform
were investigated by Opdam \cite{O} and revisited by Cherednik
\cite{C}. Here are two main results
\begin{itemize}
\item[(i)] Paley-Wiener theorem: $\kH$ and $\kI$ are (up to
positive constants) inverse isomorphisms between $C_c^\infty(\ka)$
and $P W(\kh)$.
\item[(ii)] Plancherel type formula:
$$\int_{\ka} f(x)g(-x)d\mu(x)=\text{const}\cdot \int_{i\ka} \kH f(\la)
\kH g(\la) d\nu(\la).$$
\end{itemize}
Opdam \cite{O} proved eventually a more precise Paley-Wiener
theorem: $\kH$ and $\kI$ map $C_\Gamma^\infty(\ka)$ and $P
W_{\Gamma} (\kh)$ into each other (and hence are inverse maps, up
to a positive constant), where $\Gamma$ is the convex hull of any
$W$-orbit $W\cdot x$ in $\ka$. The proof works as well for the
polar sets
$$\Gamma=\{ x\in \ka \mid (\Lambda^+,x^+)\le 1\}$$
where $\Lambda$ is any regular element in $\ka$. We shall need
this version of the Paley-Wiener theorem with positive multiples
of $\rho$.
\newline
We are now able to resume Anker's approach $\cite{A2}$ in order to
analyze the hypergeometric Fourier transform in the Schwartz
class. The following type of result was already obtained by
Delorme \cite{D}, following Harish-Chandra's strategy. On one
hand, Delorme considers only $W$-invariant functions but, on the
other hand, he deals with the more difficult case where $k<0$.
\begin{theo}
\label{inversion} The hypergeometric Fourier transform $\kH$ and
its inverse $\kI$ are topological isomorphisms between $\kC(\ka)$
and $\kS (i\ka)$.
\end{theo}
\textbf{Sketch of the proof:} The proof is divided in two parts
which correspond to the following two lemmas. The first one is
elementary.
\begin{lem}
The hypergeometric Fourier transform $\kH$ maps $\kC(\ka)$
continuously into $\kS(i\ka)$.
\end{lem}

\begin{lem}
The inverse transform $\kI: P W(i\ka) \rightarrow C_c^\infty(\ka)$
is continuous for the topology inherited from $\kS(i\ka)$ and
$\kC(\ka)$ respectively.
\end{lem}
\textbf{Proof of the lemma:} Let $h\in PW$ and $f=\kI(h)$. Given
a semi-norm $\sigma=\sigma_{p,N}$ on $\kC(\ka)$, we must find a
semi-norm $\tau$ on $\kS (i\ka)$ such that
$$\sigma_{p,N}(f) \leq \tau(h).$$
We denote by $g$ the image of $h$ by the inverse Euclidean Fourier
transform $\kF^{-1}$. According to the Paley-Wiener theorems for
the hypergeometric and the Euclidean Fourier transforms, we have
the following support conservation property: $\text{supp}(f)$ is
contained in $\Gamma_r=\{x\in \ka \mid (\rho,x^+)\le r\}$ if and
only if $\text{supp}(g) \subset \Gamma_r$. Let $\omega_j \in
C^\infty(\ka)$ such that $\omega_j=0$ inside $\Gamma_{j-1}$,
$\omega_j=1$ outside $\Gamma_j$, and $\omega_j$ is uniformly
bounded in $j\in \N^*$, as well as each derivative. Set
$g_j=\omega_jg$, $h_j=\kF(g_j)$ and $f_j=\kI(h_j)$. Here is  a
crucial observation: we have $g_j=g$ outside $\Gamma_j$, hence
$f_j=f$ outside $\Gamma_j$, by the above support property. Let us
estimate $f=f_j$ on $\Gamma_{j+1}\smallsetminus \Gamma_j$. First
of all, using Proposition \ref{derivx}, there exist $N'\in \N$ and
$C>0$ such that
$$\sup_{x\in \Gamma_{j+1}\smallsetminus \Gamma_j}(1+|x|)^NF_0(x)^{-1}|p(\frac{\partial}{\partial x})f_j(x)|<Cj^N\tau_{1,N'}(h_j).$$
Next, by the Euclidean Fourier analysis
$$\tau_{1,N'}(h_j)\leq C\sum_{l=0}^{N'}\sup_{x \in \ka}(|x|+1)^{n +1}|\nabla^l g_j(x)|.$$
Observe that $g_j$ and its derivatives vanish in $\Gamma_{j-1}$.
Hence
$$j^N\tau_{1,N'}(h_j)\leq C\sum_{l=0}^{N'}\sup_{x \in \ka \smallsetminus \Gamma_{j-1}}(|x|+1)^{N+n+1}|\nabla^l g(x)|.$$
Again, by Euclidean Fourier analysis,
$$\sum_{l=0}^{N'}\sup_{x \in \ka}(|x|+1)^{N+n+1}|\nabla^l g(x)|\leq
C\tau_{N+n+1,N''}(h).$$ In summary, there exist $N''\in \N$ and
$C>0$ such that, for every $j\in \N^*$,
\begin{eqnarray*}
\sup_{x\in \Gamma_{j+1}\smallsetminus
\Gamma_j}(1+|x|)^NF_0(x)^{-1}|p(\frac{\partial}{\partial
x})f(x)|\le C \tau_{N+n+1,N''}(h).
\end{eqnarray*}
The remaining estimate of $f$ in $\Gamma_1$ is elementary. \hfill
$\blacksquare$ \newline
\newline
In the $W$-invariant setting, the hypergeometric Fourier transform
and its inverse write $$\kH(f)(\la)= \int_\ka
f(x)F_\la(-x)d\mu(x)$$ and
$$\kI(h)(\la)=\int_{i\ka}h(\la)F_\la(x)d\nu'(\la)$$
where
\begin{eqnarray*}
d\nu'(\la)&=& \text{const}\cdot \prod_{\alpha \in \kR^+}
\frac{\Gamma((\la,\alpha^\vee)+k_\alpha+\frac{1}{2}k_{\frac{\alpha}{2}})\Gamma(-(\la,\alpha^\vee)+k_\alpha+\frac{1}{2}k_{\frac{\alpha}{2}})}{\Gamma((\la,\alpha^\vee)+\frac{1}{2}k_{\frac{\alpha}{2}})\Gamma(-(\la,\alpha^\vee)+\frac{1}{2}
k_{\frac{\alpha}{2}})}d\la \\
    &=& \text{const}\cdot \mathbf{c}(\la)^{-1}\mathbf{c}(-\la)^{-1}\ d\la
\end{eqnarray*}
is the symmetric Plancherel measure or Harish-Chandra measure (see
\cite{C2}). We denote by $\kC(\ka)^W$ and $\kS(i\ka)^W$ the spaces
of $W$-invariant functions of $\kC(\ka)$ and $\kS(i\ka)$
respectively, which we identify also with their restriction to
$\overline{\ka_+}$. From Theorem \ref{inversion} we get
\begin{cor}
\label{invsym} These transforms are topological isomorphisms
between $\kC(\ka)^W$ and $\kS(i\ka)^W$.
\end{cor}
We recover this way the main result of \cite{D} in the easy case
$k>0$.

\section{The heat kernel}

\subsection{Solution to the Cauchy problem}

In this section we solve the heat equation (with Cauchy data) for
the Heckman-Opdam Laplacian. We follow essentially the presentation of Rösler
\cite{Ros} section $4$, and refer to this article for some proofs,
which are identical in our setting. We denote by  $\kD$ the
modified Laplacian defined by
$$\kD=\frac{1}{2}(\kL-|\rho|^2).$$
The heat operator $H$ is defined by
$$H=\partial_t-\kD$$
on $C^{2,1}(\ka \times \R)$. We consider the standard Cauchy
problem: Given a continuous bounded function $f$ on $\ka$, find
$u\in C^{2,1}(\ka \times (0,+\infty))\cap C^0(\ka \times
[0,+\infty))$, such that
\begin{eqnarray} \label{cauchy} \left\{
\begin{array}{cl} H u=0 & \text{on }
 \ka \times (0,+\infty)\\
u(\cdot,0)=f.\\
\end{array}
\right. \end{eqnarray}

\begin{defin}
The heat kernel $p_t(x,y)$ is defined for $x,y\in \ka$ and $t>0$
by
\begin{eqnarray}
\label{formulenoyauchaleur}
p_t(x,y)=\int_{i\ka}e^{-\frac{t}{2}(|\la|^2+|\rho|^2)}G_\la(x)G_\la(-y)d\nu(\la).
\end{eqnarray}
The heat semigroup $(P_t,t\ge 0)$ is defined for $f\in \kC(\ka)$
and $t\ge 0$ by
$$P_tf(x):=\left\{ \begin{array}{cl}
\int_\ka p_t(x,y)f(y)d\mu(y) & \text{if }t>0 \\
f(x) & \text{if }t=0.
\end{array}
\right.$$
\end{defin}
Using the hypergeometric Fourier transform and its inverse, we can
express the heat semigroup as follows
$$P_tf=\kI(\la \mapsto
e^{-\frac{t}{2}(|\la|^2+|\rho|^2)}\kH(f)(\la))$$ and deduce its
basic properties which are summarized in the following theorem
(the analogue of Theorem $4.7$ in \cite{Ros}).
\begin{theo}
\label{semigroupe}
\begin{enumerate}
\item $(P_t,t\ge 0)$ is a strongly continuous semigroup on
$\kC(\ka)$.
\item Let $f\in \kC(\ka)$. Then $u(x,t)= P_tf(x)$
solves the Cauchy problem (\ref{cauchy}).
\end{enumerate}
\end{theo}
As in the Dunkl setting, we show next that $(P_t,t\ge 0)$ can be
extended to a strongly continuous semigroup on $C_0(\ka)$ (the
space of continuous functions $f:\ka\to \C$ which vanish at
infinity, equipped with the norm
$|f|_\infty=\sup_{x\in\ka}|f(x)|$). Consider $\kD$ as a densely
defined linear operator on $C_0(\ka)$ with domain $\kC(\ka)$.

\begin{prop}
\label{closable}
\begin{enumerate}
\item The operator $(\kD,\kC(\ka))$ has a closure, which generates a Feller semigroup $(T(t),t\ge 0)$ on
$C_0(\ka)$.
\item $T(t)$ coincides with $P_t$ on $\kC(\ka)$.
\end{enumerate}
\end{prop}
\textbf{Proof of the Proposition:} \begin{enumerate} \item In
order to apply the Hille-Yosida Theorem (see \cite{EK} Theorem
$2.2$ p.$165$) we need to check the following two properties:
\begin{enumerate}
\item Let $f \in \kC(\ka)$. Assume that $x_0$ is a global maximum of $f$. Then
$\kD f(x_0)\le 0$ (this is the positive maximum principle).
\item $(\mu I-\kD)(\kC(\ka))$ is dense in $C_0(\ka)$ for some $\mu >0$.
\end{enumerate}
$(a)$ follows from the explicit expression (\ref{explicit}) of
$\kL$. For $(b)$ we prove with Theorem \ref{inversion} that $(\mu
I-\kD)$ maps $\kC(\ka)$ onto itself for every $\mu>0$. In fact if
$f\in \kC(\ka)$, then
$$\kH((\mu I -\kD)f)(\la)=(\mu
+\frac{|\rho|^2+|\la|^2}{2})\kH(f)(\la),\ \la \in i\ka.$$
\item The equality $T(t)f=P_tf$ results from the
uniqueness of solution to (\ref{cauchy}) within the class of all
differentiable functions on $[0,\infty)$ with values in $C_0(\ka)$
(see \cite{Ros}). \hfill $\blacksquare$
\end{enumerate}
\begin{cor}
The heat kernel $p_t(x,y)$ is positive on $\ka \times \ka \times
(0,\infty)$, symmetric in $(x,y)$, and satisfies the following
properties:
\begin{enumerate}
\item For all $x,y \in \ka$, for all $t>0$ and all $w\in W$, $p_t(w x,w y)=p_t(x,y)$.
\item For all $t>0$ and $x\in \ka$, $p_t(x,\cdot) \in \kC(\ka)$.
\item Let $f\in C_b(\ka)$. Then
\begin{eqnarray*}
u(x,t)=P_tf(x)=\left\{ \begin{array}{lc} \int_\ka
p_t(x,y)f(y)d\mu(y) & \text{if } t>0 \\
f(x) & \text{if }t=0
\end{array}
\right.
\end{eqnarray*}
is still a solution to the Cauchy problem (\ref{cauchy}).
\item For all $t>0$ and all $x\in \ka$, $\int_\ka p_t(x,y)d\mu(y)=1.$
\end{enumerate}
\end{cor}
\textbf{Proof of the corollary:} The positivity property results
from the last proposition, which implies that $P_tf \ge 0$ for any
$f\in \kC(\ka)$ with $f\ge 0$. Thus (see \cite{Ros}) $p_t(x,y)\ge
0$ for all $t>0$ and $x,y\in \ka$, by continuity of
$p_t(x,\cdot)$. The invariance of $p_t$ under the Weyl group
results from the invariance of $\kD$ when $\kR^+$ is replaced by
$w\kR^+$, for any $w\in W$. The symmetry of $p_t$ results in the same way from its
invariance by $-Id$, and from Formula (\ref{formulenoyauchaleur}).
The second and third assumptions are classical and result from
basic properties of the $G_\la$ (see \cite{Ros}). The last
assumption results from the point $3$ and the fact that $T(t)1=1$
(because $\kD$ is conservative, see \cite{EK} p.$166$). \hfill
$\blacksquare$
\newline
\newline
The $W$-invariant heat kernel $p_t^W$ is defined for all $x,y\in
\ka$ and $t>0$ by
\begin{eqnarray*}
p_t^W(x,y)& =& \sum_{w\in W}p_t(x,w y)=\frac{1}{|W|}\sum_{w,w'\in W}p_t(w x,w' y)\\
          &=&
          \int_{i\ka}e^{-\frac{t}{2}(|\la|^2+|\rho|^2)}F_\la(x)F_\la(-y)d\nu'(\la).
\end{eqnarray*}
The $W$-invariant semigroup $(P_t^W,t\ge 0)$ is defined for $f \in
\kC(\overline{\ka_+})$, $x\in \overline{\ka_+}$ and $t\ge 0$, by
$$P_t^Wf(x)=\int_{\ka_+}p_t^W(x,y)f(y)d\mu(y),\ \text{if }t>0,$$
and $P^W_0f(x)=f(x)$. We have naturally the analogue of Theorem
\ref{semigroupe}. The generator of $(P_t^W,t\ge 0)$ is equal on
$\kC(\ka)^W$ to the differential part $D$ of $\kD$. The analogue
of Proposition \ref{closable} for $D$, is a consequence of
Corollary \ref{invsym} and the following lemma. The second claim
of this lemma will be used in \cite{Schap}.
\begin{lem}
The space $\kC(\ka)^W$ is dense in $C_0(\overline{\ka_+})$.
Moreover if $f\in C_c^\infty(\overline{\ka_+})$, there exists a
sequence $(u_j)_j \in \kC(\overline{\ka_+})^W$ which converges
uniformly to $f$, and which satisfies: there exists a
positive constant $C>0$, independent of $j$, such that $|\nabla u_j(x)|\le C$
for all $x\in \overline{\ka_+}$, and if
$d(x,\partial \ka_+)> \frac{1}{j}$, then $|\Delta u_j(x)|\le C$, whereas if $d(x,\partial \ka_+)\le
\frac{1}{j}$, then $|\frac{\Delta
u_j(x)}{j}|\le C$.
\end{lem}
\textbf{Proof of the lemma:} The density result is a consequence
of the Stone-Weierstrass theorem. However here we need more
information, so we need the usual technique of regularization by
convolution with an approximate of unity. Let $f\in
C_c^\infty(\overline{\ka_+})$. We can extend it to $\ka$ by
$W$-symmetry, and we get a function $\tilde{f}$ which is
symmetric, and Lipschitz. Let $u$ be an approximate of unity,
which is a $W$-symmetric $C^\infty$ function, with compact support
in the unit ball, and with integral equal to one. Then we consider
the sequence of functions $(u_j)_j$ defined by $u_j(x):=\int_\ka
\tilde{f}(x-y)j^nu(j y)dy$ for $x\in \ka$. It is classical to see
that $u_j$ is $C^\infty$ and converges uniformly to $\tilde{f}$.
It is also immediate that $u_j$ is $W$-symmetric. To see that it
has the required properties, observe that $\tilde{f}$ is derivable
(in the sense of distributions) with a bounded derivative near the
walls, and it is $C^\infty$ away from the walls.   \hfill
$\blacksquare$
\newline
\newline
We set $h_t(x)=p_t(0,x)=\frac{1}{|W|}p_t^W(0,x)$ for $x\in
\overline{\ka_+}$, and $t>0$. We have the formula:
\begin{eqnarray}
\label{noyauinvzero}
h_t(x)=\int_{i\ka}e^{-\frac{t}{2}(|\la|^2+|\rho|^2)}F_\la(x)d\nu'(\la).
\end{eqnarray}
We will now prove that the heat kernel is in fact strictly
positive. As usually we will prove this fact by using a strong
minimum principle. The result may be found in \cite{ProtWein}, but
stated in a slightly different way. Thus we include a proof.
\begin{lem}[Strong minimum principle]
\label{principemin} Let $t_0\in \R$. Let $u\in C^{2,1}(\ka \times
(t_0,+\infty))\cap C(\ka \times [t_0,+\infty))$. Assume that
$Hu\ge 0$ on $\ka\times (t_0,+\infty)$, $u\ge 0$ on $\ka \times
[t_0,+\infty)$, and $u(0,t)>0$, for all $t\ge t_0$. Then $u>0$ on
$\ka \times (t_0,+\infty)$.
\end{lem}
\textbf{Proof of the lemma:} Consider the ellipsoid
$$E\ :\ |x|^2+\gamma (t-t_0)^2< \delta.$$ Assume that $u>0$ on
$E$, and that $u(x_*,t_*)=0$ for some $(x_*,t_*)\in \partial E$,
with $t*>t_0$. By hypothesis $(x_*,t_*)$ can not be the north
pole. Moreover by reducing $E$ if necessary, we can always suppose
that it is the only point in $\overline{E}\cap\{t>t_0\}$ where $u$
vanishes. We shall perturb $u$ in a small ball
$$B\ :\ |x-x_*|^2+(t-t_*)^2<\epsilon^2,$$
with $0<\epsilon<\min(\frac{1}{2}|x_*|,\frac{1}{2}(t^*-t_0)^2)$.
Consider the auxiliary function
$$\omega(x,t)=e^{-r\delta}-e^{-r\{|x|^2+\gamma(t-t_0)^2\}}.$$
Let us compute and estimate
\begin{eqnarray*}
H\omega(x,t)&=& 2r\Big\{2r|x|^2-1+\gamma(t-t_0)-\sum_{\alpha \in
\kR^+}k_\alpha (\alpha,x)\coth(\frac{\alpha}{2},x)\Big\}\\
            &\times& e^{-r\{|x|^2+\gamma(t-t_0)^2\}}.
\end{eqnarray*}
This expression can be made strictly positive on $\overline{B}$,
by choosing $r>0$ sufficiently large. The function
$v=u+\epsilon'\omega$
\begin{itemize}
\item[$\bullet$] is strictly positive on $\overline{B}\setminus
\overline{E}$, since $\omega>0$ outside of $\overline{E}$,
\item[$\bullet$] is equal to $u$ on $\overline{B}\cap\partial E$,
since $\omega$ vanishes on $\partial E$,
\item[$\bullet$] can be made strictly positive on $\partial B \cap
\overline{E}$ by choosing $\epsilon'>0$ sufficiently small.
\end{itemize}
Thus the minimum $v_*\le 0$ of $v$ on $\overline{B}$ is achieved
at an inner point. There $\partial_tv=0$, $\nabla v=0$, and
$\Delta v\ge 0$. Hence $Hv\le0$. But on the other side
$Hv=Hu+\epsilon'H\omega>0$, and we have a contradiction. \hfill
$\blacksquare$
\newline
\newline
We can deduce from this lemma the
\begin{cor}
The heat kernel $p_t(x,y)$ is strictly positive on $\ka \times \ka
\times (0,+\infty)$.
\end{cor}
\textbf{Proof of the corollary:} First we apply the preceding
lemma for the function $u(x,t)=h_t(x)$. We have simply to prove
that $h_t(0)$ is strictly positive for all $t>0$. This comes from
Formula (\ref{noyauinvzero}). Moreover since the preceding lemma
may be applied for any $t_0>0$, we get that $p_t(x,0)>0$ for any
$t>0$ and $x\in \ka$. Suppose now that $p_t(x,y)=0$ for some
$x,y\in \ka-\{0\}$ and $t>0$. We have
$$p_t(x,y)=\int_\ka
p_{\frac{t}{2}}(x,z)p_{\frac{t}{2}}(z,y)d\mu(z).$$ But as $p$ is
positive and continuous, this implies that
$p_{\frac{t}{2}}(x,0)p_{\frac{t}{2}}(0,y)=0$, and we get a
contradiction. \hfill $\blacksquare$
\begin{rem} \emph{
Since the space $\kC(\ka)$ is dense in all the $L^p(\ka,\mu)$
spaces, for $p\in [1,\infty)$, the Hille-Yosida theorem (cf
\cite{EK}) assures that $\kD$ is closable on $L^p(\ka,\mu)$ and
generates a heat semigroup $(T^{(p)}(t),t\ge 0)$, which is
strongly continuous. Moreover, still by an argument of uniqueness
in the Cauchy problem, we see that $T^{(p)}$ coincides with the
preceding operator $P$ on $\kC(\ka)$. And by continuity we see
that $T^{(p)}$ is just the natural extension of $P$ on
$L^p(\ka,\mu)$. It is equal for $f\in L^p(\ka,\mu)$, $x\in \ka$,
and $t>0$, to
$$T^{(p)}(t)f(x)=P_tf(x)=\int_\ka p_t(x,y)f(y)d\mu(y).$$
Obviously the same discussion apply in the radial situation (with
$D$ and $P^W$ in place of $\kD$ and $P$ respectively). }
\end{rem}

\subsection{Estimates and asymptotic of the heat kernel}
In this subsection we establish a sharp global estimate of $h_t$
(Theorem \ref{estimchaleur}) and an asymptotic of
$p_T(x,\sqrt{T}y)$ when $T\to \infty$ (Proposition
\ref{asymptotiquechaleur2}). Let $\gamma:=\sum_{\alpha \in
\kR_+}k_\alpha$, and like usually for $x\in \ka$, we denote by
$x^+$ its symmetric in $\overline{\ka_+}$. A problem in order to
get global estimate of $p_t$ is that it is not a convolution
operator. Thus $p_t(\cdot,\cdot)$ can not be simply expressed in
terms of the function $h_t(\cdot)$. Therefore the next Theorem is
only a partial result. A better one could be obtain if we had a
global estimate of the Dunkl kernel.
\begin{theo}
\label{estimchaleur} The following global estimate holds, for all
$t>0$ and $x\in \ka$:
\begin{eqnarray*} h_t(x) &\asymp& t^{-\gamma-\frac{n}{2}}\left\{
\prod_{\alpha \in
\kR_0^+}(1+|(\alpha,x)|)(1+t+|(\alpha,x)|)^{k_\alpha+k_{2\alpha}-1}
\right\} \\
 & \times & e^{-|\rho|^2\frac{t}{2}-(\rho,x^+)-\frac{|x|^2}{2t}}.
\end{eqnarray*}
\end{theo}
\textbf{Proof of the theorem:} Thanks to Theorem \ref{euler}, and
the known expression of the heat kernel associated to the Dunkl
Laplacian \cite{Ros}, we can use exactly the same proof as in
\cite{AO}. In this proof it was made use of the heat kernel in
balls of radius $R>0$ with boundary conditions. This may be
avoided by using weak parabolic minimum (or maximum) principles
for unbounded domains, which hold also because the heat kernel
vanishes at infinity. \hfill $\blacksquare$
\newline
\newline
Our next result gives an equivalent of $p_T(x,\sqrt{T}y)$ when $T$
tends to $\infty$. This result will be needed in \cite{Schap} for
the proof of the convergence of the normalized $F_0$-process.
However since the proof is easier in the $W$-invariant case, we
begin by the analogue result for $p_T^W(x,\sqrt{T}y)$. Then we
will simply explain what has to be modified in the non invariant
setting.
\begin{prop}
\label{asymptotiquechaleur} There exists a constant $K>0$, such
that for any $x \in \overline{\ka_+}$ and any $y\in \ka_+$,
$$p_T^W(x,\sqrt{T}y) \sim K
e^{-\frac{|y|^2}{2}}T^{-\frac{n}{2}-|\kR_0^+|}e^{-\frac{|\rho|^2}{2}T}F_0(-x)F_0(\sqrt{T}y),$$
when $T\to +\infty$.
\end{prop}
\textbf{Proof of the proposition:} We resume the "analysis away
from walls" carried out in \cite{AJ}. It consists in expanding
$F_{\la}$ in the heat kernel expression
\begin{eqnarray}
\label{intsemi}
p_T^W(x,\sqrt{T}y)=\int_{i\ka}e^{-\frac{T}{2}(|\la|^2+|\rho|^2)}F_\la(-x)F_\la(\sqrt{T}y)d\nu'(\la)
\end{eqnarray}
using the Harish-Chandra series \cite{HO}
$$F_\la(y)=\sum_{w\in W}\mathbf{c}(w\la)e^{(w\la-\rho,y)}\sum_{q \in Q^+}\Gamma_q(w\la)
e^{-(q,y)}.$$ Recall that this expression holds for $y\in \ka_+$.
Now we replace $F_\la(\sqrt{T}y)$ by its development in series in
the integral (\ref{intsemi}). The properties of the coefficients
$q_\chi$ allow us to invert the integral term and the series (see
\cite{AJ} for more details). Therefore we get
\begin{eqnarray}
\label{seriesemi}p_T^W(x,\sqrt{T}y)=\sum_{q \in Q^+}
E_q(x,y)e^{-\frac{|\rho|^2}{2}T-(\rho+q,\sqrt{T}y)}
\end{eqnarray}
where (using the $W$-invariance of $\nu'$ in $\la$), for $x,y\in
\ka$,
\begin{eqnarray*}
E_q(x,y) = K \int_{i\ka}e^{-\frac{T}{2}|\la|^2+(\la,\sqrt{T}y)}
F_\la(-x)c(\la)\Gamma_q(\la)d\nu'(\la).
\end{eqnarray*}
Here $K$ is a constant whose value may change in the sequel. We
denote by  $\mathbf{b'}$ the function defined by
$$\mathbf{b'}(\la)\frac{\mathbf{c}(\la)}{\mathbf{\pi}(\la)}d\nu'(\la)=d\la.$$ It is holomorphic
in zero. Observe now that
$$\mathbf{\pi}(\frac{1}{T}\frac{\partial}{\partial_\la})e^{-\frac{T}{2}|\la|^2}=\mathbf{\pi}(-\la)e^{-\frac{T}{2}|\la|^2}.$$
This formula comes from the fact that there are no skew symmetric polynomial of strictly lower degree than
$|\kR_0^+|$. Thus the function $E_0$ may be rewritten into
$$E_0(x,y)=
K\int_{i\ka}e^{-\frac{T}{2}|\la|^2}\pi(\frac{1}{T}\frac{\partial}{\partial_\la})\{e^{(\la,\sqrt{T}y)}
F_\la(-x)b'(\la)^{-1}\}d\la.$$ Then we make the change of
variables $v=\frac{y+\la}{\sqrt{T}}$, and we get
\begin{eqnarray*}
E_0(x,y)e^{-\frac{|\rho|^2}{2}T-(\rho,\sqrt{T}y)} &\sim &  K
F_0(-x)e^{-\frac{T}{2}|\rho|^2-(\rho,\sqrt{T}y)-\frac{|y|^2}{2}}T^{-\frac{D}{2}}\pi(\sqrt{T}y)\\
&\times & \int_{i\ka}e^{\frac{1}{2}|v|^2}
\frac{F_{\frac{v-y}{\sqrt{T}}}(-x)}{F_0(-x)}{b'}^{-1}(\frac{v-y}{\sqrt{T}})dv.
\end{eqnarray*}
The preceding integral has a finite limit, independent of $x$ and
$y$, when $T$ tends to infinity. Thus using the known asymptotic
of $F_0$ (Theorem \ref{estimate1}), we conclude that the first
term of the series in (\ref{seriesemi}) has the desired
asymptotic. A similar study would show that the leading terms are
negligible. This concludes the proof of the proposition. \hfill
$\blacksquare$
\begin{prop}
\label{asymptotiquechaleur2} There exists a constant $K>0$, such
that for any $x \in \ka$, and any $y\in \ka_{\text{reg}}$, if $wy\in
\ka_+$, then
$$p_T(x,\sqrt{T}y) \sim K
e^{-\frac{|y|^2}{2}}T^{-\frac{n}{2}-|\kR_0^+|}e^{-\frac{|\rho|^2}{2}T}G_0(w
x)F_0(\sqrt{T}y),$$ when $T\to +\infty$.
\end{prop}
\textbf{Proof of the proposition:} The proof is analogue as for
the preceding proposition. First we have
$p_T(x,\sqrt{T}y)=p_T(wx,w\sqrt{T}y)$. Then in the integral
expression of $p_T(wx,w\sqrt{T}y)$ , we replace $G_\la (-w
\sqrt{T} y)$ by its development in series. Since $-wy\in \ka_-$,
we already know the dominant coefficients of the development. Indeed they
were computed by Opdam in \cite{O}: they are all null except one
which is equal up to a constant to $\pi(\la)$. But
$\pi(\la)d\nu(\la)$ behaves like $d\nu'(\la)$ in zero, i.e. like
$|\pi(\la)|^2$. Thus we can follow the rest of the proof of the
preceding proposition, and we get the result. \hfill
$\blacksquare$

\subsection{The Poisson equation for $\kD$}
Our sharp estimates of Theorem \ref{estimderivees} allows us
to prove the
\begin{prop} Let $f\in L^1(\ka,\mu)$. Then the function
$Gf:x\mapsto \int_0^\infty P_tf(x)dt$ is finite $\mu$-a.e. If
moreover $\kF(f) \in L^1(i\ka,\nu)$, then $Gf$ is bounded, belongs
to $C^2(\ka)$, and satisfies the Poisson equation $\kD Gf=-f$.
\end{prop}
\textbf{Proof of the proposition:} Let $f\in L^1(\ka,\mu)$. For
all $x$, and all $\epsilon
>0$, we have
\begin{eqnarray*} |G f(x)|&=& |\int_0^\infty
e^{-\frac{t}{2}|\rho|^2}\int_\ka \int_{i\ka}
e^{-\frac{t}{2}|\la|^2}G(\la,x)G(-\la,y)d\nu(\la)f(y)d\mu(y)dt|\\
                         &\le & |\int_0^1 P_tf(x)dt|+C|f|_1\int_1^\infty
                         e^{-\frac{t}{2}|\rho|^2}
                         |\int_{i\ka}e^{-\frac{t}{2}|\la|^2}d\nu(\la)|dt,
\end{eqnarray*}
where $C$ is a constant. But since for any $t\ge 0$, $P_t$ is a
contraction on $L^1$, we have therefore $|P_tf|_1\le |f|_1$. Thus
$|\int_0^1 P_tf dt|_1 \le |f|_1<\infty$. And then $\mu$-a.e.,
$|\int_0^1 P_tfdt|<\infty$. Finally we get that $\mu$ a.e.
$Gf<\infty$. This proves the first claim of the proposition. Now let $f\in L^1(\ka,\mu)$, be such that $\kF(f) \in L^1(i\ka,\nu)$.
Then we have
\begin{eqnarray*}
|G f(x)| &\le & \int_{i\ka}\kF(f)(\la)\int_0^\infty
e^{-\frac{t}{2}|(|\la|^2+\rho|^2)}dt d\nu(\la)\\
&\le & 2\int_{i\ka}\frac{\kF(f)(\la)}{|\la|^2+|\rho|^2}d\nu(\la).
\end{eqnarray*}
This shows that $Gf$ is bounded. Moreover using a theorem of derivation
under the integral, and our precise estimate of the derivatives
of the functions $G_\la$, we see that $Gf \in C^2(\ka)$ and
satisfies $\kD Gf=-f$. This finishes the proof of the proposition.
\hfill $\blacksquare$

\section{Appendix : computation of the Heckman-Opdam laplacian}

First we give another expression of the Cherednik operator:
\begin{eqnarray*}
T_\xi f(x)&=&\partial_\xi f(x)+\sum_{\alpha\in
\kR^+}\frac{k_\alpha}{2}(\alpha,\xi)\coth\frac{(\alpha,x)}{2}\{f(x)-f(r_\alpha
.x)\}\\
&-& \sum_{\alpha \in
\kR^+}\frac{k_\alpha}{2}(\alpha,\xi)f(r_\alpha.x).
\end{eqnarray*}
Now we compute
\begin{eqnarray*}
& &\sum_{\alpha\in \kR^+}
\frac{k_\alpha}{2}\coth\frac{(\alpha,x)}{2}\{T_\alpha
f(x)-T_\alpha f(r_\alpha
.x)\}\\
&=& \sum_{\alpha\in
\kR^+}\frac{k_\alpha}{2}\coth\frac{(\alpha,x)}{2}\{\partial_\alpha
f(x)-\partial_\alpha f(r_\alpha
.x)\}\\
&+&\sum_{\alpha,\beta \in \kR^+} \frac{k_\alpha
k_\beta}{4}(\alpha,\beta)\coth\frac{(\alpha,x)}{2}\coth\frac{(\beta,x)}{2}\{f(x)-f(r_\alpha
.x)\}\\
&-&\sum_{\alpha,\underbrace{\beta}_{\beta'}\in
\kR^+}\frac{k_\alpha
\overbrace{k_\beta}^{k_\beta'}}{4}\overbrace{(\alpha,\beta)}^{-(\alpha,\beta')}
\coth\frac{(\alpha,x)}{2}\coth\frac{\overbrace{(\beta,r_\alpha .
x)}^{(\beta',x)}}{2}\{f(r_\alpha .x)-f(\overbrace{r_\beta
r_\alpha}^{r_\alpha r_{\beta'}}
.x)\}\\
&-&\sum_{\alpha,\beta \in \kR^+} \frac{k_\alpha
k_\beta}{4}(\alpha,\beta)\coth\frac{(\alpha,x)}{2}f(r_\beta .x)\\
&+& \sum_{\alpha,\beta \in \kR^+} \frac{k_\alpha
k_\beta}{4}(\alpha,\beta)\coth\frac{(\alpha,x)}{2}f(r_\beta
r_\alpha .x)\\
&=&\sum_{\alpha\in
\kR^+}\frac{k_\alpha}{2}\coth\frac{(\alpha,x)}{2}\{\partial_\alpha
f(x)-\partial_\alpha f(r_\alpha
.x)\}\\
&+&\sum_{\alpha,\beta \in \kR^+} \frac{k_\alpha
k_\beta}{4}(\alpha,\beta)\coth\frac{(\alpha,x)}{2}\coth\frac{(\beta,x)}{2}\{f(x)-f(r_\beta
r_\alpha
.x)\}\\
&-&\sum_{\alpha,\beta \in \kR^+} \frac{k_\alpha
k_\beta}{4}(\alpha,\beta)\coth\frac{(\alpha,x)}{2}\{f(r_\beta .x)
-f(r_\beta r_\alpha .x)\}.
\end{eqnarray*}
Thanks to the following lemma, we can remove the hyperbolic
cotangent in the second sum.
\begin{lem}
Let $\kR$ be an integral root system (non necessarily reduced).
Then
$$\sum_{\alpha,\beta \in \kR^+,r_\beta \circ
r_\alpha=\tau}k_\alpha k_\beta
(\alpha,\beta)\{\coth\frac{(\alpha,x)}{2}\coth\frac{(\beta,x)}{2}-1\}=0$$
for all non trivial rotation $\tau$.
\end{lem}
\textbf{Proof of the lemma:} Applying the Euclidean Laplacian to
the Weyl denominator formula
$$\prod_{\alpha \in
\kR^+}\{e^{\frac{(\alpha,x)}{2}}-e^{-\frac{(\alpha,x)}{2}}\}=\sum_{w\in
W}e^{(w. \rho,x)}$$ we get the identity
$$\sum_{\alpha,\beta \in \kR^+,
\alpha\neq\beta}(\alpha,\beta)\{\coth\frac{(\alpha,x)}{2}\coth\frac{(\beta,x)}{2}-1\}=0$$
which holds for all reduced root system. Now by restricting to the
different root systems of rank $2$, we see that this relation is
equivalent to the lemma.
\newline
$\bullet$ \textbf{$\mathbf{A_1 \times A_1}$:} trivial.
\newline
$\bullet$ \textbf{$\mathbf{A_2}$:}
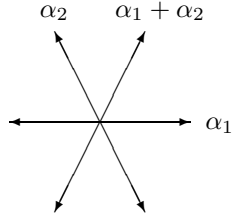
\begin{center}
\begin{picture}(8,7) \put(4,3){\vector(1,0){3}}
\put(4,3){\vector(1,2){1.5}} \put(4,3){\vector(-1,2){1.5}}
\put(4,3){\vector(-1,0){3}} \put(4,3){\vector(-1,-2){1.5}}
\put(4,3){\vector(1,-2){1.5}} \put(7.5,2.75){$\alpha_1$}
\put(2,6.5){$\alpha_2$} \put(4.5,6.5){$\alpha_1+\alpha_2$}
\end{picture}
\end{center}
The lemma reduces to the identity
$$\frac{k^2}{2}\{-\coth\frac{\alpha_1}{2}\coth\frac{\alpha_2}{2}+\coth\frac{\alpha_1}{2}\coth\frac{\alpha_1+\alpha_2}{2}+\coth\frac{\alpha_2}{2}\coth\frac{\alpha_1+\alpha_2}{2}-1\}=0.$$
$\bullet$ \textbf{$\mathbf{B_2=C_2}$:}
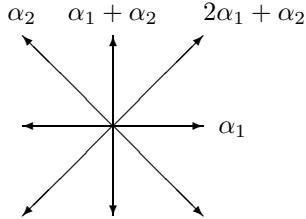
\begin{center}
\begin{picture}(8,7) \put(4,3){\vector(1,0){3}}
\put(4,3){\vector(1,1){3}} \put(4,3){\vector(-1,1){3}}
\put(4,3){\vector(0,1){3}} \put(4,3){\vector(-1,0){3}}
\put(4,3){\vector(-1,-1){3}} \put(4,3){\vector(1,-1){3}}
\put(4,3){\vector(0,-1){3}} \put(7.5,2.75){$\alpha_1$}
\put(0.5,6.5){$\alpha_2$} \put(2.5,6.5){$\alpha_1+\alpha_2$}
\put(7,6.5){$2\alpha_1+\alpha_2$}
\end{picture}
\end{center}
The lemma reduces to the identity
\begin{eqnarray*}
k_1k_2\{&-&
\coth\frac{\alpha_1}{2}\coth\frac{\alpha_2}{2}+\coth\frac{\alpha_1}{2}\coth\frac{2\alpha_1+\alpha_2}{2}\\
&+&
\coth\frac{\alpha_2}{2}\coth\frac{\alpha_1+\alpha_2}{2}+\coth\frac{\alpha_1+\alpha_2}{2}\coth\frac{2\alpha_1+\alpha_2}{2}-2\}=0.
\end{eqnarray*}
$\bullet$ \textbf{$\mathbf{BC_2}$:}
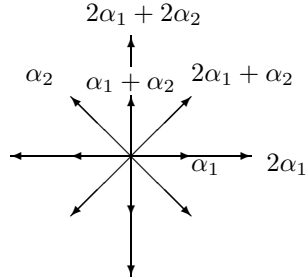
\begin{center}
\begin{picture}(10,9) \put(5,4){\vector(1,0){2}}
\put(5,4){\vector(1,1){2}} \put(5,4){\vector(-1,1){2}}
\put(5,4){\vector(0,1){2}} \put(5,4){\vector(-1,0){2}}
\put(5,4){\vector(-1,-1){2}} \put(5,4){\vector(1,-1){2}}
\put(5,4){\vector(0,-1){2}} \put(7,4){\vector(1,0){2}}
\put(3,4){\vector(-1,0){2}}\put(5,2){\vector(0,-1){2}}\put(5,7){\vector(0,1){1}}
\put(7,3.5){$\alpha_1$} \put(9.5,3.5){$2\alpha_1$}
\put(1.5,6.5){$\alpha_2$} \put(3.5,6.25){$\alpha_1+\alpha_2$}
\put(7,6.5){$2\alpha_1+\alpha_2$}\put(3.5,8.5){$2\alpha_1+2\alpha_2$}
\end{picture}
\end{center}
The lemma reduces to the following identities of type $B_2=C_2$
\begin{eqnarray*}
k_1k_2\{&-&
\coth\frac{\alpha_1}{2}\coth\frac{\alpha_2}{2}+\coth\frac{\alpha_1}{2}\coth\frac{2\alpha_1+\alpha_2}{2}\\
&+&
\coth\frac{\alpha_2}{2}\coth\frac{\alpha_1+\alpha_2}{2}+\coth\frac{\alpha_1+\alpha_2}{2}\coth\frac{2\alpha_1+\alpha_2}{2}-2\}=0\\
2k_2k_3\{&-&
\coth\frac{\alpha_1}{2}\coth\frac{\alpha_2}{2}+\coth\frac{\alpha_1}{2}\coth\frac{2\alpha_1+\alpha_2}{2}\\
&+&
\coth\frac{\alpha_2}{2}\coth\frac{\alpha_1+\alpha_2}{2}+\coth\frac{\alpha_1+\alpha_2}{2}\coth\frac{2\alpha_1+\alpha_2}{2}-2\}=0.
\end{eqnarray*}
$\bullet$ \textbf{$\mathbf{G_2}$:}
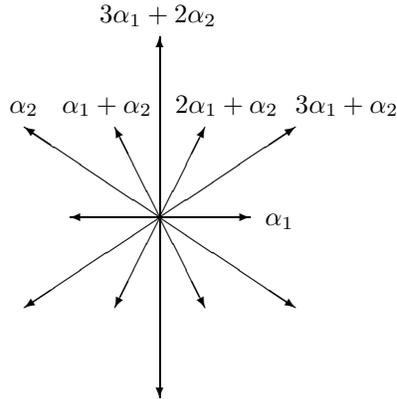
\begin{center}
\begin{picture}(12,11) \put(5,5){\vector(1,0){3}}
\put(5,5){\vector(1,2){1.5}} \put(5,5){\vector(-1,2){1.5}}
\put(5,5){\vector(-1,0){3}} \put(5,5){\vector(-1,-2){1.5}}
\put(5,5){\vector(1,-2){1.5}} \put(5,5){\vector(3,2){4.5}}
\put(5,5){\vector(-3,-2){4.5}} \put(5,5){\vector(0,1){6}}
\put(5,5){\vector(0,-1){6}} \put(5,5){\vector(-3,2){4.5}}
\put(5,5){\vector(3,-2){4.5}} \put(8.5,4.75){$\alpha_1$}
\put(0,8.5){$\alpha_2$} \put(1.75,8.5){$\alpha_1+\alpha_2$}
\put(5.5,8.5){$2\alpha_1+\alpha_2$}
\put(9.5,8.5){$3\alpha_1+\alpha_2$}
\put(3,11.5){$3\alpha_1+2\alpha_2$}
\end{picture}
\end{center}
The lemma reduces to the following identities, the last ones being
of type $A_2$
\begin{eqnarray*}
\frac{3k_1k_2}{2}\{&-&
\coth\frac{\alpha_1}{2}\coth\frac{\alpha_2}{2}+\coth\frac{\alpha_1}{2}\coth\frac{3\alpha_1+\alpha_2}{2}\\
&+&
\coth\frac{\alpha_2}{2}\coth\frac{\alpha_1+\alpha_2}{2}+\coth\frac{\alpha_1+\alpha_2}{2}\coth\frac{3\alpha_1+2\alpha_2}{2}\\
&+&\coth\frac{2\alpha_1+\alpha_2}{2}\coth\frac{3\alpha_1+\alpha_2}{2}+\coth\frac{2\alpha_1+\alpha_2}{2}\coth\frac{3\alpha_1+2\alpha_2}{2}
-4\}=0\\
\frac{k_1^2}{2}\{&-&
\coth\frac{\alpha_1}{2}\coth\frac{\alpha_1+\alpha_2}{2}+\coth\frac{\alpha_1}{2}\coth\frac{2\alpha_1+\alpha_2}{2}\\
&+&
\coth\frac{\alpha_1+\alpha_2}{2}\coth\frac{2\alpha_1+\alpha_2}{2}-1\}=0\\
\frac{3k_2^2}{2}\{&-&\coth\frac{\alpha_2}{2}\coth\frac{3\alpha_1+\alpha_2}{2}+\coth\frac{\alpha_2}{2}\coth\frac{3\alpha_1+2\alpha_2}{2}\\
&+&
\coth\frac{3\alpha_1+\alpha_2}{2}\coth\frac{3\alpha_1+2\alpha_2}{2}-1\}=0.
\end{eqnarray*}
\hfill $\blacksquare$
\newline
\newline
Eventually we get the expression of the Heckman-Opdam Laplacian:
\begin{eqnarray*}
\kL f(x)&=&\sum_{j=1}^n T_{\xi_j}^2f(x)\\
&=&\sum_{j=1}^n \partial_{\xi_j}T_{\xi_j}f(x)+ \sum_{\alpha \in
\kR^+}\frac{k_\alpha}{2}\coth\frac{(\alpha,x)}{2}\overbrace{\sum_{j=1}^n(\alpha,\xi_j)\{T_{\xi_j}f(x)-T_{\xi_j}f(r_\alpha
.x)\}}^{T_\alpha f(x)-T_\alpha f(r_\alpha .x)}\\
&-& \sum_{\alpha \in
\kR^+}\frac{k_\alpha}{2}\underbrace{\sum_{j=1}^n(\alpha,\xi_j)T_{\xi_j}f(r_\alpha
.x)}_{T_\alpha f(r_\alpha .x)} \\
&=&\Delta f(x) + \sum_{\beta \in
\kR^+}\frac{k_\beta}{4}\underbrace{\sum_{j=1}^n(\beta,\xi_j)^2}_{|\beta|^2}\overbrace{(1-\coth^2\frac{(\beta,x)}{2})}^{-\sinh^{-2}\frac{(\beta,x)}{2}}\{f(x)-f(r_\beta
.x)\} \\
&+&\sum_{\beta \in \kR^+}\frac{k_\beta}{2}\coth
\frac{(\beta,x)}{2}\overbrace{\sum_{j=1}^n
(\beta,\xi_j)\{\partial_{\xi_j}f(x)-\partial_{r_\beta .
\xi_j}f(r_\beta .x)\}}^{\partial_\beta f(x)+\partial_\beta
f(r_\beta .x)}\\
&-&\sum_{\beta \in
\kR^+}\frac{k_\beta}{2}\underbrace{\sum_{j=1}^n(\beta,\xi_j)\partial_{r_\beta
.\xi_j}f(r_\beta .x)}_{-\partial_\beta f(r_\beta .x)}\\
&+&\sum_{\alpha \in \kR^+}\frac{k_\alpha}{2}\coth
\frac{(\alpha,x)}{2}\{\partial_\alpha f(x)-\partial_\alpha
f(r_\alpha .x)\}\\
&+&\sum_{\alpha,\beta \in \kR^+}\frac{k_\alpha
k_\beta}{4}(\alpha,\beta)\{f(x)-f(r_\beta r_\alpha .x)\}\\
&-&\sum_{\alpha,\beta \in \kR^+}\frac{k_\alpha
k_\beta}{4}(\alpha,\beta)\coth\frac{(\alpha,x)}{2}\{f(r_\beta .x)-f(r_\beta r_\alpha .x)\}\\
&-&\sum_{\alpha,\underbrace{\beta}_{\beta'} \in
\kR^+}\frac{k_\alpha
\overbrace{k_\beta}^{k_\beta'}}{4}\overbrace{(\alpha,\beta)}^{-(\alpha,\beta')}\coth\frac{\overbrace{\beta,r_\alpha
.x)}^{(\beta',x)}}{2}\{f(r_\alpha .x)-f(\overbrace{r_\beta
r_\alpha}^{r_\alpha r_{\beta'}} .x))\}\\
&+&\sum_{\alpha,\beta \in \kR^+} \frac{k_\alpha
k_\beta}{4}(\alpha,\beta)f(r_\beta r_\alpha .x)\\
&=& \Delta f(x) + \sum_{\alpha \in \kR^+}k_\alpha \coth
\frac{(\alpha,x)}{2}\partial_\alpha f(x) + |\rho|^2f(x)\\
&-&\sum_{\alpha \in \kR^+}k_\alpha
\frac{|\alpha|^2}{4\sinh^2\frac{(\alpha,x)}{2}}\{f(x)-f(r_\alpha
.x)\}.
\end{eqnarray*}

\end{document}